\documentclass{amsart}

\usepackage{amssymb}
\usepackage{amsmath}
\usepackage{amsthm}
\usepackage{amsfonts}

\usepackage{a4wide}

\usepackage{lineno}

\newtheorem{lemma}{Lemma}
\newtheorem{proposition}{Proposition}
\newtheorem{theorem}{Theorem}
\newtheorem{corollary}{Corollary}

\theoremstyle{definition}
\newtheorem{definition}{Definition}
\newtheorem{example}{Example}

\newtheorem{remark}{Remark}
\newtheorem{conjecture}{Conjecture}

\numberwithin{equation}{section}

\begin{document}


\title[Algebraic number fields generated by monogenic trinomials]
{Algebraic number fields generated by \\ an infinite family of monogenic trinomials}

\author{Daniel C. Mayer}
\address{Naglergasse 53\\8010 Graz\\Austria}
\email{algebraic.number.theory@algebra.at}
\urladdr{http://www.algebra.at}

\author{Abderazak Soullami}
\address{Department of Mathematics, Faculty of Sciences and Technology\\Moulay Ismail University, Errachidia\\Morocco}
\email{abderazak.soullami@usmba.ac.ma}

\dedicatory{Dedicated to the memory of Georgi F. Voronoi}

\thanks{Research of the first author supported by the Austrian Science Fund (FWF): projects P26008-N25 and J0497-PHY,
and by the Research Executive Agency of the European Union (EUREA)}



\date{Sunday, 10 April 2022}


\begin{abstract}
For an infinite family of monogenic trinomials
\(P(X)=X^3\pm 3rbX-b\in\mathbb{Z}\lbrack X\rbrack\),
arithmetical invariants of the cubic number field \(L=\mathbb{Q}(\theta)\),
generated by a zero \(\theta\) of \(P(X)\),
and of its Galois closure \(N=L(\sqrt{d_L})\) are determined.
The conductor \(f\) of the cyclic cubic relative extension \(N/K\),
where \(K=\mathbb{Q}(\sqrt{d_L})\) denotes the unique quadratic subfield of \(N\),
is proved to be of the form \(3^eb\) with \(e\in\lbrace 1,2\rbrace\),
which admits statements concerning
primitive ambiguous principal ideals, lattice minima,
and independent units in \(L\).
The number \(m\)
of non-isomorphic cubic fields \(L_1,\ldots,L_m\)
sharing a common discriminant \(d_{L_i}=d_L\)
with \(L\) is determined.
\end{abstract}

\maketitle


\section{Introduction}
\label{s:Intro}

\noindent
Parametrized families of monogenic polynomials have been given in 
\cite{BSS2022,SmSl2020}.
In the Sections \S\S\
\ref{s:Conductors},
\ref{s:Multiplets}
of the present article,
we investigate the arithmetic of the cubic number field
\(L=\mathbb{Q}(\theta)\)
generated by a zero \(\theta\) of a polynomial
\begin{equation}
\label{eqn:Trinomial}
P(X)=X^{3}+aX-b, \quad \text{ with } a=\sigma\cdot 3rb, \quad \text{ sign } \sigma\in\lbrace -1,+1\rbrace, \quad \text{ and } r,b\in\mathbb{N},
\end{equation}
which arises by specialization from the trinomials
\(X^{p^n}+aX^{p^s}-b\), \(p\in\mathbb{P}\), \(n>s\ge 0\), in
\cite{BSS2022}.
The cubic field \(L\) is \textit{simply real} for \(\sigma=+1\) and \textit{totally real} for \(\sigma=-1\).

Many of the number theoretic properties of the cubic field \(L\)
and its maximal order \(\mathcal{O}_L\)
are governed by the \textit{relative conductor} \(f=f_{N/K}\)
of the Galois closure \(N=L(\sqrt{d_L})\) of \(L\),
viewed as a cyclic cubic extension of the unique quadratic subfield
\(K=\mathbb{Q}(\sqrt{d_L})\)
\cite{Ha1930}.
When the trinomial \(P(X)\) is \textit{monogenic}, we prove that
\begin{equation}
\label{eqn:Conductor}
f=3^e\cdot b, \quad \text{ with exponent } e\in\lbrace 1,2\rbrace, 
\end{equation}
whence the ramification of primes in \(L\)
can be tuned intentionally by selection of the polynomial coefficient \(b\).
Thus, the trinomials in Equation
\eqref{eqn:Trinomial}
are very useful for applications in Section \S\ 
\ref{s:Application}.

In particular,
since the norm of a primitive ambiguous principal ideal of \(N/K\)
is a divisor of the square \(f^2\),
the \textit{differential principal factorization} (DPF) in \(L\) and \(N\)
\cite{Ma2019,Ma2021a,Ma2021b}
is closely connected with the coefficient \(b\).
When \(b\) is prime or composite,
we show that the zero \(\theta\) of \(P(X)\) is generator of
the ideal \(\mathfrak{Q}\) of \(\mathcal{O}_L\)
whose third power \(\mathfrak{Q}^3\) is equal to \(b\mathcal{O}_L\). Thus,
\begin{equation}
\label{eqn:DPF}
\mathfrak{Q}=b\mathcal{O}_L+\theta\mathcal{O}_L=\theta\mathcal{O}_L
\end{equation}
is an ambiguous principal ideal and constitutes an \textit{absolute} DPF in \(L\),
which forces \(L\) to be of DPF-\textit{type} \(\beta\) if \(\sigma=+1\),
and of one of the DPF-\textit{types} \(\beta_1,\beta_2,\gamma,\varepsilon\) if \(\sigma=-1\)
\cite{Ma2021b,Mo1979}.
In the simply real situation,
\(\theta\), respectively \(1-3r\theta=\theta^3/N_{L/\mathbb{Q}}(\theta)\),
is a \textit{fundamental unit} of \(L\), if \(b=1\), respectively \(b>1\).
In infinitely many cases of \(r\), \(\theta\) and \(\theta^2\) are
\textit{lattice minima} of \(\mathcal{O}_L\)
\cite{Vo1896},
for each fixed \(b>1\).

By means of the ray class field routines by Fieker
\cite{Fi2001},
which are implemented in the computational algebra system Magma
\cite{BCP1997,BCFS2022,MAGMA2022},
it is possible to construct the complete multiplet \((L_1,\ldots,L_m)\)
of all non-isomorphic cubic fields sharing the common discriminant
\(d_{L_i}=d_L\), for all \(1\le i\le m\)
\cite{Ha1930,Ma1992,Ma2014}.
If \(u\) free and \(v\) restrictive primes divide \(f\), and \(\omega=v_3(b)\),
the multiplicity
\begin{equation}
\label{eqn:Multiplicity}
m=3^{\varrho+\omega}\cdot 2^{u}\cdot\frac{1}{3}\lbrack 2^{v-1}-(-1)^{v-1}\rbrack, \quad
\text{ where } \varrho\in\mathbb{N}_0 \text{ denotes the \(3\)-class rank of } K,
\end{equation}
is in intimate relationship to the number \(\tau=u+v\) of prime divisors of \(f\),
growing exponentially with the latter.
We show that the monogenic field \(L\)
with square free \(b\)
is a distinguished member of its multiplet \((L_1,\ldots,L_m)\),
because it is unique with mandatory \textit{absolute principal factorization}.

Before we actually restrict ourselves to the cubic situation \(p=3\)
we give general criteria for generators of ambiguous principals in Section \S\
\ref{s:Ambiguous}.
For numerical verifications see Sections \S\S\
\ref{s:Experiments},
\ref{s:ExperimentsReal}.


\section{Ambiguous principal ideals}
\label{s:Ambiguous}

\noindent
Let \(E/F\) be an arbitrary (not necessarily Galois) extension of algebraic number fields
with relative degree \(\lbrack E:F\rbrack=d\). Denote by
\(\iota:\,\mathcal{I}_F\to\mathcal{I}_E\), \(\mathfrak{a}\mapsto\mathfrak{a}\mathcal{O}_E\),
the embedding of the ideal group of the base field \(F\)
into the ideal group of the extension field \(E\) with maximal order \(\mathcal{O}_E\).

An ideal \(\mathfrak{A}\in\mathcal{I}_E\) of the extension field \(E\)
is called \textit{ambiguous} with respect to \(E/F\),
if its \(d\)-th power \(\mathfrak{A}^d\in\iota(\mathcal{I}_F)\)
is an embedded ideal of the base field \(F\).


\begin{lemma}
\label{lem:PowerByNorm}
For an algebraic number \(A\in E^\times\) of the extension field \(E\) and
the principal ideal \(A\mathcal{O}_E\in\mathcal{P}_E\) generated by \(A\),
the following necessary and sufficient criterion holds.
\begin{equation}
\label{eqn:Ambiguity}
(\exists\,\mathfrak{n}\in\mathcal{I}_F)\,(A\mathcal{O}_E)^d=\iota(\mathfrak{n}) \quad\Longleftrightarrow\quad
(\exists\,U\in U_E)\,\frac{A^d}{\mathrm{N}_{E/F}(A)}=U,
\end{equation}
that is, the principal ideal \(A\mathcal{O}_E\) is ambiguous with respect to \(E/F\)
if and only if the quotient of the power \(A^d\) by the norm \(\mathrm{N}_{E/F}(A)\) is a unit of \(E\).
\end{lemma}

\begin{proof}
If the quotient \textit{power by norm}, \(\frac{A^d}{\mathrm{N}_{E/F}(A)}=U\in U_E\),
of the algebraic integer \(A\) is a unit in \(E\),
then the power \(A^d=\nu U\) belongs to the group \(F^\times\cdot U_E\),
where \(\nu:=\mathrm{N}_{E/F}(A)\in F^\times\).
Consequently, the principal ideal
\((A\mathcal{O}_E)^d=A^d\mathcal{O}_E=\nu U\mathcal{O}_E=\nu\mathcal{O}_E=\iota(\nu\mathcal{O}_F)\)
is contained in the embedding \(\iota(\mathcal{P}_F)<\iota(\mathcal{I}_F)\),
whence the principal ideal \(A\mathcal{O}_E\) is ambiguous with respect to \(E/F\).
Conversely, the \textit{ambiguity} of the principal ideal \(A\mathcal{O}_E\)
implies the existence of an ideal \(\mathfrak{n}\in\mathcal{I}_F\) such that 
\((A\mathcal{O}_E)^d=\iota(\mathfrak{n})=\mathfrak{n}\mathcal{O}_E\).
Taking the norm, we get the power relation
\(\mathrm{N}_{E/F}(A\mathcal{O}_E)^d=\mathrm{N}_{E/F}((A\mathcal{O}_E)^d)
=\mathrm{N}_{E/F}(\mathfrak{n}\mathcal{O}_E)=\mathfrak{n}^d\)
in the ideal group \(\mathcal{I}_F\simeq\mathbb{Z}^{(\mathbb{P}_F)}\) of \(F\)
which is torsion free abelian, whence \(\mathrm{N}_{E/F}(A\mathcal{O}_E)=\mathfrak{n}\).
Thus, the ideal \(\mathfrak{n}\) is principal, generated by the norm \(\mathrm{N}_{E/F}(A)\).
Finally, the equality of principal ideals \(A^d\mathcal{O}_E=(A\mathcal{O}_E)^d=\mathrm{N}_{E/F}(A)\mathcal{O}_E\)
enforces the existence of a unit \(U\in U_E\) such that
\(A^d=\mathrm{N}_{E/F}(A)\cdot U\).
Thus the quotient power by norm, \(\frac{A^d}{\mathrm{N}_{E/F}(A)}=U\), is a unit.
\end{proof}


\noindent
Let \(F\) be an algebraic number field
with maximal order \(\mathcal{O}_F\).

\begin{theorem}
\label{thm:TrinomialPrinciple}
Let \(P(X)=X^d+aX^t-b\in\mathcal{O}_F\lbrack X\rbrack\)
be a monic irreducible trinomial with exponents \(d>t\ge 1\)
and middle coefficient \(a=vb\) divisible by the absolute coefficient \(b\).
In the extension field \(E=F(\theta)\)
of degree \(\lbrack E:F\rbrack=d\) of \(F\),
generated by a zero \(\theta\) of \(P(X)\),
\begin{equation}
\label{eqn:TrinomialPrinciple}
\begin{aligned}
& \theta \text{ is generator of an \textbf{ambiguous principal ideal} with norm } \mathrm{N}_{E/F}(\theta)=(-1)^{d-1}b, \text{ and} \\
& U=\frac{\theta^d}{\mathrm{N}_{E/F}(\theta)}=(-1)^{d-1}(1-v\theta^t)\in U_E
\text{ is a \textbf{unit}}.
\end{aligned}
\end{equation}
\end{theorem}

\begin{proof}
The monic irreducible trinomial \(P(X)=X^d+aX^t-b\) with algebraic integer coefficients in \(F\),
and exponents \(d>t\ge 1\) is supposed to have
a middle coefficient \(a=vb\) divisible by the absolute coefficient \(b\).
If the conjugates of the zero \(\theta\) of \(P(X)\) are denoted by
\(\theta=\theta^{(0)},\theta^{(1)},\theta^{(2)},\ldots,\theta^{(d-1)}\),
\begin{equation*}
\begin{aligned}
\text{then } \qquad
P(X) &= (X-\theta^{(0)})(X-\theta^{(1)})(X-\theta^{(2)})\cdots(X-\theta^{(d-1)}) \\
     &= X^d+aX^t+(-1)^{d}\theta^{(0)}\theta^{(1)}\theta^{(2)}\cdots\theta^{(d-1)} \\
     &= X^d+aX^t+(-1)^{d}\mathrm{N}_{E/F}(\theta).
\end{aligned}
\end{equation*}
So the absolute coefficient \(-b=(-1)^{d}\mathrm{N}_{E/F}(\theta)\)
is essentially the norm of the zero \(\theta\), up to the sign.
Consequently, the quotient \textit{power by norm} of the zero \(\theta\) (with sign \(+\) for odd \(d\))
is given by
\[
\frac{\theta^d}{\mathrm{N}_{E/F}(\theta)}=\frac{b-a\theta^t}{\pm b}=\frac{b-vb\theta^t}{\pm b}
=\pm(1-v\theta^t)\in\mathcal{O}_F\lbrack\theta\rbrack.
\]
Since the relative equation order \(\mathcal{O}_F\lbrack\theta\rbrack\) is contained in
the maximal order \(\mathcal{O}_E\) of the extension \(E\),
the number \(\pm(1-v\theta^t)\) is an algebraic integer with norm
\[
\mathrm{N}_{E/F}\left((-1)^{d-1}(1-v\theta^t)\right)
=\mathrm{N}_{E/F}\left(\frac{\theta^d}{\mathrm{N}_{E/F}(\theta)}\right)
=\frac{\mathrm{N}_{E/F}(\theta^d)}{\mathrm{N}_{E/F}(\mathrm{N}_{E/F}(\theta))}
=\frac{\mathrm{N}_{E/F}(\theta)^d}{\mathrm{N}_{E/F}(\theta)^d}=1.
\]
Thus, \(U=(-1)^{d-1}(1-v\theta^t)\) is a unit of \(E\),
and \(\theta\) generates an ambiguous principal ideal of \(E/F\), by Lemma
\ref{lem:PowerByNorm}.
\end{proof}


\begin{remark}
\label{rmk:TrinomialPrinciple}
If the absolute coefficient \(b\) in Theorem
\ref{thm:TrinomialPrinciple}
is distinct from \(\pm 1\),
then the zero \(\theta\) is generator of a \textit{non-trivial} ambiguous principal ideal.
However, if \(b=\pm 1\), then \(\theta\in U_E\) itself is a unit of \(E\),
and \(U=(-1)^{d-1}(1-v\theta^t)=\pm\theta^d\)
is certainly \textit{not} an \textit{independent} unit
(in particular, \textit{not} a \textit{fundamental} unit in the case of unit rank \(1\)).

Another special case arises for \(v=0\) and thus \(a=0\).
Then we have a \textit{pure} equation \(P(X)=X^d-b=0\) and
\(\theta=\sqrt[d]{b}\) is a \textit{radical}.
In this case, the unit \(U=(-1)^{d-1}=\pm 1\) is \textit{trivial}.

Up to now, nothing was required concerning the prime decomposition of 
the absolute coefficient \(b\) in Theorem
\ref{thm:TrinomialPrinciple}.
It was neither assumed to be squarefree nor to have bounded exponents.
However, in order to get a \textit{primitive} ambiguous principal ideal
\(\theta\mathcal{O}_E\)
generated by the zero \(\theta\), it should be assumed that
\(b=\pm\mathrm{N}_{E/F}(\theta)\) is \textit{free of \(d\)-th powers}.
\end{remark}


\noindent
Now we draw a conclusion about the trinomials
given by Boughaleb, Soullami, Sahmoudi in
\cite{BSS2022},
\textit{without} imposing the condition of \textit{monogeneity}.

\begin{corollary}
\label{cor:BSS}
Let \(F\) be an algebraic number field
with maximal order \(\mathcal{O}_F\), and
let \(p\in\mathbb{P}\) be a rational prime number.
Suppose that \(P(X)=X^{p^n}+aX^{p^s}-b\in\mathcal{O}_F\lbrack X\rbrack\)
is a monic irreducible but \textbf{not necessarily monogenic}
trinomial with \(n>s\ge 0\)
and middle coefficient \(a=vb\) divisible by the absolute coefficient \(b\).
In the extension field \(E=F(\theta)\)
of degree \(\lbrack E:F\rbrack=p^n\) of \(F\),
generated by a zero \(\theta\) of \(P(X)\),
\begin{equation}
\begin{aligned}
\label{eqn:BSS}
& \theta \text{ is generator of an \textbf{ambiguous principal ideal} with norm } \mathrm{N}_{E/F}(\theta)=(-1)^{p^n-1}b, \text{ and} \\
& U=\frac{\theta^{p^n}}{\mathrm{N}_{E/F}(\theta)}=(-1)^{p^n-1}(1-v\theta^{p^s})\in U_E
\text{ is a \textbf{unit}}.
\end{aligned}
\end{equation}
\end{corollary}

\begin{proof}
The claims are an immediate consequence of Theorem
\ref{thm:TrinomialPrinciple}
by specialization of the exponents
\(d=p^n\) and \(t=p^s\) with \(n>s\ge 0\)
and thus \(d>t\ge 1\), as required.
In order that Formula
\eqref{eqn:TrinomialPrinciple}
can be applied,
we have to assume the divisibility  \(a=vb\) again.
\end{proof}


\begin{remark}
\label{rmk:BSS}
Further specialization to interesting situations within \textit{experimental} reach
by computational algebra systems like Magma
\cite{MAGMA2022,BCP1997,BCFS2022}
is possible for the rational base field \(F=\mathbb{Q}\),
\(a=vb\),
and some \textit{composite exponents} such as
\begin{enumerate}
\item
two types of \textit{nonic} extensions with \(P(X)=X^9+aX^{3^s}-b\),
\(p=3\), \(n=2\) and \(s\in\lbrace 0,1\rbrace\),
\item
three types of \textit{octic} extensions with \(P(X)=X^8+aX^{2^s}-b\),
\(p=2\), \(n=3\) and \(s\in\lbrace 0,1,2\rbrace\),
\item
two types of \textit{quartic} extensions with \(P(X)=X^4+aX^{2^s}-b\),
\(p=2\), \(n=2\) and \(s\in\lbrace 0,1\rbrace\),
\end{enumerate}
and especially for \textit{prime exponents} such as investigated in
\cite[Section \S\ 5]{BSS2022},
where \(v\) is assumed to be of the particular form \(v=\pm pr\), i.e.,
\(P(X)=X^p\pm prbX-b\)
with \(n=1\) and \(s=0\).
Since it will turn out that the Galois group of \(P(X)\)
is the \textit{full symmetric group} \(S_p\) for \(p\ge 5\),
we shall restrict to the cubic situation \(p=3\),
where complete theories of ambiguous principal ideals
and multiplets sharing a common conductor are available,
in the Sections \S\S\
\ref{s:Conductors},
\ref{s:Multiplets}.
\end{remark}


\section{Discriminants and conductors of cubic number fields}
\label{s:Conductors}

\begin{theorem}
\label{thm:Preparation}
\textbf{(Preparation Theorem.)}
\begin{itemize}
\item
For a trinomial \(P(X)=X^3+aX-b\)
with \(a=\sigma\cdot 3rb\), \(\sigma\in\lbrace -1,+1\rbrace\), \(r,b\in\mathbb{N}\),
the cubic number field \(L=\mathbb{Q}(\theta)\),
generated by a zero \(P(\theta)=0\) of \(P(X)\),
possesses the \textbf{discriminant}
\[
d_L=d_P/i(\theta)^2, \text{ where } d_P=-27b^2\partial \text{ with } \partial=\sigma\cdot 4r^3b+1,
\]
and \(i(\theta)=(\mathcal{O}_L:\mathbb{Z}\lbrack\theta\rbrack)\)
denotes the index of the equation order in the maximal order of \(L\).
The field \(L\) is complex (simply real) if \(\sigma=+1\),
and \(L\) is totally real (triply real) if \(\sigma=-1\).
\item
If \(P(X)\) is \textbf{monogenic} with square free coefficient \(b\),
then \(i(\theta)=1\) and \(d_L=d_P=-27b^2\partial\),
and for the decomposition
\[d_L=f^2\cdot d_K\]
of \(d_L\) into the square of the \textbf{conductor} \(f\)
of the Galois closure \(N=L(\sqrt{d_L})\) of \(L\)
over the unique quadratic subfield \(K=\mathbb{Q}(\sqrt{d_L})\)
and the fundamental discriminant \(d_K\) of \(K\)
there exist three situations:
\begin{equation}
\label{eqn:Main}
\begin{aligned}
f=3b,\ b=3f_0,\ m=2,\quad 3\nmid\partial,\ d_K=-3\partial,      \qquad &\text{ if } 3\parallel b \textbf{ (irregular case),} \\
f=9b,\ b= f_0,\ m=2,\quad 3\parallel\partial,\ d_K=-\partial/3, \qquad &\text{ if } 3\nmid b,\ 3\mid(\sigma\cdot rb+1), \\
f=3b,\ b= f_0,\ m=1,\quad 3\nmid\partial,\ d_K=-3\partial,      \qquad &\text{ if } 3\nmid b,\ 3\nmid(\sigma\cdot rb+1),
\end{aligned}
\end{equation}
where \(f_0\) denotes the squarefree part of \(f\) coprime to \(3\)
and \(f=3^m\cdot f_0\) with \(m\in\lbrace 1,2\rbrace\).
In each situation,
\(\partial\) is the \textbf{dual discriminant} of the dual quadratic field
\(\mathbb{Q}(\sqrt{\partial})\) of \(K\),
in the sense of the \textbf{cubic reflection theorem} by Scholz
\cite{So1932}.
\end{itemize}
\end{theorem}

\begin{proof}
\(\bullet\)
By specialization of
\cite[Lemma 2 and Corollary 3]{BSS2022},
the discriminant of the polynomial \(P(X)=X^3+aX-b\),
that is \(X^{p^n}+aX^{p^s}-b\) with \(p=3\), \(n=1>s=0\),
is given by
\(d_P=(-1)^{\frac{p^n(p^n-1)}{2}}\cdot p^{sp^n}b^{p^s-1}\Delta^{p^s}=-\Delta\)
where the \textit{core discriminant} is defined by
\(\Delta=p^pb^{p-1}+(p-1)^{p-1}a^p=27b^2+4a^3\).
In the particular situation with
\(a=\sigma\cdot 3rb\), \(\sigma\in\lbrace -1,+1\rbrace\), \(r,b\in\mathbb{N}\),
we obtain
\(\Delta=27b^2+4\cdot 27\sigma^3r^3b^3=27b^2\partial\)
with \textit{dual quadratic discriminant} \(\partial=1+\sigma\cdot 4r^3b\),
since the sign satisfies the relation \(\sigma^3=\sigma\).
Obviously, \(\partial>0\) if \(\sigma=+1\), and \(\partial<0\) if \(\sigma=-1\).
From \(a=\sigma\cdot 3rb\), it follows that
\(v_3(a)\ge v_3(b)+1\) and \(v_p(a)\ge v_p(b)\),
for \(p\in\mathbb{P}\setminus\lbrace 3\rbrace\).

\(\bullet\)
According to
\cite[Theorem 1]{BSS2022},
\(P(X)\) is \textit{monogenic} if and only if \(b\) is square free and
\(v_3(T)=1\) for the \textit{critical term} \(T=b^2+a-1=b^2+\sigma\cdot 3rb-1\).
The essential part of the proof is now conducted with the aid of
\cite[Thm. 2, p. 583]{LlNt1983} by Llorente and Nart,
which also depends on the same \textit{critical term} \(T\).
Note that \(v_3(d_L)=2\) cannot occur,
since \(3\nmid d_K\) and \(3\mid f\) \(\Longrightarrow\) \(v_3(d_L)=4\).
\begin{enumerate}
\item
\(v_3(d_L)=5\) if and only if \(1=v_3(b)<v_3(b)+1\le v_3(a)\), i.e. \(3\mid b\).
\item
\(v_3(d_L)=4\) if and only if \(3\mid a\), \(v_3(a+3)\ge 2\), \(3\nmid b\), \(v_3(T)=1\).
\item
\(v_3(d_L)=3\) if and only if \(3\mid a\), \(v_3(a+3)=1\), \(3\nmid b\), \(v_3(T)=1\).
\item
\(v_3(d_L)=1\) is excluded, since neither \(v_3(b)>1\) nor \(v_3(T)\ge 2\).
\item
\(v_3(d_L)=0\) is impossible, since neither \(3\nmid a\) nor \(v_3(T)\ge 3\).
\end{enumerate}
According to
\cite[Thm. 3, p. 584]{LlNt1983} by Llorente and Nart,
which has been proved earlier with different methods as
\cite[Satz 6, p. 578]{Ha1930} by Hasse,
we have \(d_L=f^2\cdot d_K\) with \(f=3^mf_0\), \(m\in\lbrace 0,1,2\rbrace\),
\(f_0\) square free and coprime to \(3\).

For primes \(p\ge 5\), we have
\(p\mid f_0\) \(\iff\) \(v_p(d_L)=2\) \(\iff\) \(1=v_p(b)\le v_p(a)\) \(\iff\) \(p\mid b\).
For the prime \(p=2\), we have
\(2\mid f_0\) \(\iff\) \(v_2(d_L)=2\) and \(2\nmid d_K\) \(\iff\) \(1=v_2(b)\le v_2(a)\) \(\iff\) \(2\mid b\).
In the case \(2\mid f_0\), we must have \(d_K\equiv 5\,(\mathrm{mod}\,8)\).
Thus, \(b=f_0\) if \(3\nmid b\) and \( b=3f_0\) if \(3\mid b\).

Now we analyze the \(3\)-valuation of the discriminant
\(d_L=3^{2m}f_0^2\cdot d_K=-27b^2\partial\), for which generally
\(v_3(d_L)=2m+v_3(d_K)=3+2v_3(b)+v_3(\partial)\).
We use \(a+3=\sigma\cdot 3rb+3=3(\sigma\cdot rb+1)\).

\begin{enumerate}
\item
If \(v_3(b)=1\), then \(5=v_3(d_L)=2m+v_3(d_K)=5+v_3(\partial)\),
which is only possible for \(m=2\), \(v_3(d_K)=1\), \(v_3(\partial)=0\), 
that is \(f=9f_0=3b\) and \(d_K=-3\partial\equiv -3\,(\mathrm{mod}\,9)\).
\item
If \(v_3(b)=0\) and \(3\mid(\sigma\cdot rb+1)\), then \(4=v_3(d_L)=2m+v_3(d_K)=3+v_3(\partial)\),
and consequently \(m=2\), \(v_3(d_K)=0\), \(v_3(\partial)=1\),
that is \(f=9f_0=9b\) and \(d_K=-\partial/3\).
\item
If \(v_3(b)=0\) and \(3\nmid(\sigma\cdot rb+1)\), then \(3=v_3(d_L)=2m+v_3(d_K)=3+v_3(\partial)\),
which implies \(m=1\), \(v_3(d_K)=1\), \(v_3(\partial)=0\),
and thus \(f=3f_0=3b\) and \(d_K=-3\partial\equiv\pm 3\,(\mathrm{mod}\,9)\).
\end{enumerate}
In each situation, \(d_K=-3\partial/\gcd(3,\partial)^2\) is the \textit{dual discriminant} of \(\partial\).
\end{proof}


\section{Principal factors and multiplicities of cubic number fields}
\label{s:Multiplets}

\noindent
In the case of monogeneity, we have the following general \textbf{principal ideal criterion}.

\begin{lemma}
\label{lem:Monogenic}
Let \(L=\mathbb{Q}(\theta)\)
be an algebraic number field of degree
\(\lbrack L:\mathbb{Q}\rbrack=n\),
generated by adjunction of a zero \(\theta\)
of a monic irreducible polynomial
\[
P(X)=X^n+c_{n-1}X^{n-1}+c_{n-2}X^{n-2}+\ldots+c_2X^2+c_1X+c_0\in\mathbb{Z}\lbrack X\rbrack
\]
to the field \(\mathbb{Q}\) of rational numbers.
Denote by \(\mathcal{O}\) the maximal order of \(L\) and suppose that
\(\mathfrak{Q}=q\mathcal{O}+\theta\mathcal{O}\)
is an ideal of \(\mathcal{O}\) with a rational integer \(q\in\mathbb{Z}\).
If \(P(X)\) is \textbf{monogenic} and
the absolute coefficient \(c_0\) of \(P(X)\) divides \(q\), \(c_0\mid q\),
then \(\mathfrak{Q}\) is the \textbf{principal ideal} \(\theta\mathcal{O}\) generated by \(\theta\).
\end{lemma}

\begin{proof}
By the assumption of \textbf{monogeneity},
the zero \(\theta\) of the monic irreducible polynomial
\(P(X)=\sum_{i=0}^{n}\,c_iX^i\) with \(c_n=1\)
is generator of a \textbf{power basis}
\((1,\theta,\theta^{2},\ldots,\theta^{n-1})\)
for the maximal order \(\mathcal{O}\) of \(L\) over the ring \(\mathbb{Z}\) of rational integers, i.e.
\[
\mathcal{O}=\mathbb{Z}\oplus\mathbb{Z}\theta\oplus\mathbb{Z}\theta^{2}\oplus\ldots\oplus\mathbb{Z}\theta^{n-1}.
\]
In the monic equation \(0=P(\theta)=\sum_{i=0}^{n}\,c_i\theta^i\),
the leading coefficient \(c_n=1\) admits a representation of the biggest power as
\(\theta^n=-\sum_{i=0}^{n-1}\,c_i\theta^i\).

Since an element \(\theta\mu\) of the principal ideal \(\theta\mathcal{O}\) can be viewed as
\(q\cdot 0+\theta\mu\), we have the inclusion
\(\theta\mathcal{O}\subset q\mathcal{O}+\theta\mathcal{O}\).
In order to prove the reverse inclusion, we proceed as follows.
An element of the ideal \(\mathfrak{Q}=q\mathcal{O}+\theta\mathcal{O}\) has the general form
\(q\varkappa+\theta\lambda\) with
\(\varkappa=\sum_{i=0}^{n-1}\,k_i\theta^i\in\mathcal{O}\) and \(\lambda=\sum_{i=0}^{n-1}\,\ell_i\theta^i\in\mathcal{O}\),
and we obtain
\begin{equation*}
\label{eqn:LHS}
\begin{aligned}
q\varkappa+\theta\lambda &= q\sum_{i=0}^{n-1}\,k_i\theta^i + \theta\sum_{i=0}^{n-1}\,\ell_i\theta^i \\
                         &= qk_0 + \sum_{i=1}^{n-1}\,qk_i\theta^i + \sum_{i=0}^{n-2}\,\ell_i\theta^{i+1} + \ell_{n-1}\theta^n \\
                         &= qk_0 + \sum_{i=1}^{n-1}\,qk_i\theta^i + \sum_{j=1}^{n-1}\,\ell_{j-1}\theta^{j} - \ell_{n-1}\sum_{i=0}^{n-1}\,c_i\theta^i \\
                         &= (qk_0-\ell_{n-1}c_0) + \sum_{j=1}^{n-1}\,(qk_j+\ell_{j-1}-\ell_{n-1}c_j)\theta^j.                      
\end{aligned}
\end{equation*}
We want to establish the existence of an algebraic integer
\(\mu=\sum_{i=0}^{n-1}\,m_i\theta^i\in\mathcal{O}\)
such that
\(q\varkappa+\theta\lambda=\theta\mu\).
Thus, we compute
\begin{equation*}
\label{eqn:RHS}
\begin{aligned}
\theta\mu &= \theta\sum_{i=0}^{n-1}\,m_i\theta^i \\
          &= \sum_{i=0}^{n-2}\,m_i\theta^{i+1} + m_{n-1}\theta^n \\
          &= \sum_{j=1}^{n-1}\,m_{j-1}\theta^{j} - m_{n-1}\sum_{i=0}^{n-1}\,c_i\theta^i \\
          &= -m_{n-1}c_0 + \sum_{j=1}^{n-1}\,(m_{j-1}-m_{n-1}c_j)\theta^j.
\end{aligned}
\end{equation*}
Comparing the coefficients with respect to the power basis yields
\begin{equation*}
\label{eqn:Coeff}
\begin{aligned}
-m_{n-1}c_0        &= qk_0-\ell_{n-1}c_0, \text{ and } \\
m_{j-1}-m_{n-1}c_j &= qk_j+\ell_{j-1}-\ell_{n-1}c_j, \text{ for all } 1\le j\le n-1, \text{ respectively} \\
m_{j-1}            &= qk_j+\ell_{j-1}+(m_{n-1}-\ell_{n-1})c_j, \text{ for all } 1\le j\le n-1.
\end{aligned}
\end{equation*}
Now we need the assumption that \(c_0\mid q\).
The solubility by means of \(\mu\in\mathcal{O}\)
depends decisively on the integral solution
\(m_{n-1}=\ell_{n-1}-\frac{q}{c_0}k_0\).
\end{proof}


\begin{remark}
\label{rmk:Monogenic}
The condition \(c_0\mid q\) in Lemma
\ref{lem:Monogenic}
is satisfied by all \textbf{monogenic} trinomials
\(P(X)=X^{p^m}+aX^{p^s}-b\) with \(p\in\mathbb{P}\), \(m>s\ge 0\)
in
\cite{BSS2022},
under the identification \(b=c_0=q\).
However, in order that
\(\mathfrak{Q}=q\mathcal{O}+\theta\mathcal{O}\)
can be interpreted as a prime ideal \(\mathfrak{Q}\in\mathbb{P}_{\mathcal{O}}\),
we must additionally assume a prime \(b=q\in\mathbb{P}\), and \(q\mid a\) is needed
to get the congruence \(P(X)\equiv X^{n}\,(\mathrm{mod}\,q)\) with \(n=p^m\).
To be independent of the requirement of a \textit{prime} coefficient \(b\),
we have stated another necessary and sufficient criterion for an
absolutely ambiguous principal ideal (absolute DPF)
in Lemma
\ref{lem:PowerByNorm}.
\end{remark}


\begin{theorem}
\label{thm:Main}
\textbf{(Main Theorem.)}
Under the assumptions of Theorem
\ref{thm:Preparation},
a zero \(\theta\) of the \textbf{monogenic} polynomial
\(P(X)=X^3+\sigma\cdot 3rbX-b\in\mathbb{Z}\lbrack X\rbrack\)
with sign \(\sigma\in\lbrace -1,+1\rbrace\),
\textbf{squarefree} coefficient \(b\in\mathbb{N}\)
and parameter \(r\in\mathbb{N}\), such that
\(v_3(b^2+\sigma\cdot 3rb-1)=1\) and
\((\forall\,\ell\in\mathbb{P})\) \(\lbrack v_\ell(\sigma\cdot 4r^3b+1)\ge 2\Longrightarrow\ell\mid 3rb\rbrack\),
generates a cubic field \(L=\mathbb{Q}(\theta)\) with maximal order
\(\mathcal{O}_L=\mathbb{Z}\lbrack\theta\rbrack\),
discriminant \(d_L=-27b^2(\sigma\cdot 4r^3b+1)\), conductor either \(f=3b\) or \(f=9b\),
and the following additional arithmetic properties.
\begin{enumerate}
\item
If \(b\) is a prime or composite number, then
\(\mathfrak{Q}=b\mathcal{O}_L+\theta\mathcal{O}_L\in\mathcal{I}_L\) is an ambiguous principal ideal.
Since \(\mathfrak{Q}\) is an \textbf{absolute differential principal factor} of \(L/\mathbb{Q}\),
the cubic field \(L\) must be of type \(\beta\) if \(L\) is simply real
and of one of the types \(\beta_1,\beta_2,\gamma,\varepsilon\) if \(L\) is triply real.
\item
If the quadratic subfield \(K\) of the Galois closure of \(L\)
has trivial \(3\)-class group \(\mathrm{Cl}_3(K)=1\),
and \(L\) is simply real,
then the \textbf{multiplicity} \(m\) of the homogeneous multiplet \((L_1,\ldots,L_m)\)
of non-isomorphic cubic fields sharing the discriminant \(d_{L_i}=d_L\), for all \(1\le i\le m\),
is given by \(m=3^\omega 2^{\tau-1}\),
where \(\tau\) denotes the number of prime divisors of the conductor \(f\),
and \(\omega=1\) if \(9\parallel f\), \(d_K\equiv -3\,(\mathrm{mod}\,9)\), but \(\omega=0\) otherwise
(in fact, \(\omega=v_3(b)\)).
\item
If \(L\) is simply real and \(b>1\), then \(\varepsilon_0=1-3r\theta\) is a \textbf{fundamental unit} of \(L\),
i.e. \(U_L=\langle -1,\varepsilon_0\rangle\),
but if \(b=1\), the zero \(\theta\) itself is a fundamental unit of \(L\), \(U_L=\langle -1,\theta\rangle\).
\end{enumerate}
\end{theorem}

\begin{proof}
According to
\cite[Lem. 1, p. 581]{LlNt1983},
the polynomial decomposition with three identical factors \(G(X)=X\) in
\(P(X)=X^3+\sigma\cdot 3rbX-b\equiv X^3=X\cdot X\cdot X\,(\mathrm{mod}\,q_i)\)
modulo each prime divisor \(q_i\in\mathbb{P}\) of \(b=q_1\cdots q_k\)
establishes a prime factorization
\(q_i\mathcal{O}_L=\mathfrak{Q}_i^3\)
into the third power of the prime ideal
\(\mathfrak{Q}_i=q_i\mathcal{O}_L+\theta\mathcal{O}_L\in\mathbb{P}_L\) of \(\mathcal{O}_L\),
where trivially \(\theta=G(\theta)\),
that is, \(q_i\) is totally ramified in \(L\).
In the special case that \(-c_0=b=q\) is itself prime, we can immediately apply Lemma
\ref{lem:Monogenic}.
Together with the fact that \(b\mid f\) is a divisor of the conductor \(f=3^eb\),
the Lemma implies that
\(\mathfrak{Q}=q\mathcal{O}_L+\theta\mathcal{O}_L=\theta\mathcal{O}_L\)
is the ambiguous principal ideal generated by \(\theta\).
If \(b=q_1\cdots q_k\) is composite with \(k\ge 2\),
then \(\mathfrak{Q}=b\mathcal{O}_L+\theta\mathcal{O}_L=\mathfrak{Q}_1\cdots\mathfrak{Q}_k\).
For \(k\ge 3\), this follows by induction from the subsequent proof for \(k=2\):
if \(\mathfrak{Q}_1=q_1\mathcal{O}_L+\theta\mathcal{O}_L\in\mathbb{P}_L\), 
\(\mathfrak{Q}_2=q_2\mathcal{O}_L+\theta\mathcal{O}_L\in\mathbb{P}_L\), and
\(\mathfrak{Q}=q_1q_2\mathcal{O}_L+\theta\mathcal{O}_L\), then
\(\mathfrak{Q}\subset\mathfrak{Q}_1\) and \(\mathfrak{Q}\subset\mathfrak{Q}_2\), i.e.,
\(\mathfrak{Q}_1\) and \(\mathfrak{Q}_2\) are two prime ideals dividing \(\mathfrak{Q}\),
and thus \(\mathfrak{Q}\subset\mathfrak{Q}_1\cap\mathfrak{Q}_2=\mathfrak{Q}_1\mathfrak{Q}_2\).
Conversely, \(\mathfrak{Q}_1\mathfrak{Q}_2\subset\mathfrak{Q}\), since
\(\xi=q_1\alpha+\theta\beta\in\mathfrak{Q}_1\) and \(\eta=q_2\gamma+\theta\delta\in\mathfrak{Q}_2\)
implies \(\xi\eta=q_1q_2\alpha\gamma+\theta(q_1\alpha\delta+q_2\beta\gamma+\theta\beta\delta)\in\mathfrak{Q}\).
With \(-c_0=b=q_1\cdots q_k\) the Lemma establishes that
\(\mathfrak{Q}=b\mathcal{O}_L+\theta\mathcal{O}_L=\theta\mathcal{O}_L\)
is the ambiguous principal ideal generated by \(\theta\).
In each case,
\(\mathfrak{Q}^3=\mathfrak{Q}_1^3\cdots\mathfrak{Q}_k^3=q_1\cdots q_k\mathcal{O}_L=b\mathcal{O}_L\),
the zero \(\theta\) is an \textbf{absolute principal factor} of \(L\)
\cite{Ma2019,Ma2021a,Ma2021b},
and \(L\) must be of type \(\beta\) if \(L\) is simply real,
and of one of the types \(\beta_1,\beta_2,\gamma,\varepsilon\) if \(L\) is totally real.
In the first case, the types \(\alpha_1,\alpha_2\),
and in the last case, the types \(\alpha_1,\alpha_2,\alpha_3,\delta_1,\delta_2\) are excluded,
since they are free of absolute principal factors.
(They only possess relative principal factors
\cite{Ma2019,Ma2021a,Ma2021b}
or capitulation
\cite{Ar1927,Ar1929}.)

According to
\cite{Ma1992,Ma2014},
the defect \(\delta=0\) vanishes and thus the multiplicity is
\(m=3^\omega2^{\tau-1}\).

It remains to show that there exists an algebraic integer
\(\varepsilon_0^{-1}=e_0+e_1\theta+e_2\theta^2\in\mathcal{O}_L\), such that
\begin{equation*}
\begin{aligned}
1 &= \varepsilon_0\cdot\varepsilon_0^{-1}=(1-3r\theta)(e_0+e_1\theta+e_2\theta^2) \\
  &= e_0+e_1\theta+e_2\theta^2-3re_0\theta-3re_1\theta^2-3re_2\theta^3 \\
  &= e_0+(e_1-3re_0)\theta+(e_2-3re_1)\theta^2-3re_2(b-a\theta) \\
  &= (e_0-3rbe_2)+(e_1+3rae_2-3re_0)\theta+(e_2-3re_1)\theta^2 \\ 
  &= (e_0-3rbe_2)+(e_1+3^2r^2be_2-3re_0)\theta+(e_2-3re_1)\theta^2 \\ 
  &= (e_0-3rbe_2)+(e_1+3r(3rbe_2-e_0))\theta+(e_2-3re_1)\theta^2.
\end{aligned}
\end{equation*}
Comparing the coefficients with respect to the power integral basis yields
\begin{equation*}
\begin{aligned}
 e_1 &= -3r(3rbe_2-e_0) = -3r\cdot (-1) = 3r, \\
 e_2 &= 3re_1 = 3r\cdot 3r = 9r^2, \\
 e_0 &= 1+3rbe_2 = 1+3rb\cdot 9r^2 = 1+27r^3b.
\end{aligned}
\end{equation*}
So the \textbf{inverse fundamental unit} is given by
\begin{equation}
\label{eqn:Inverse}
\varepsilon_0^{-1}=(1+27r^3b)+3r\theta+9r^2\theta^2.
\end{equation}

\noindent
We must calculate the minimal polynomial \(M_{\varepsilon}\) of the unit \(\varepsilon=\varepsilon_0\).
The polynomial must be monic with absolute coefficient the negative norm
\(\mathrm{N}_{L/\mathbb{Q}}(\varepsilon)=1\).
Thus, we try the approach
\begin{equation}
\label{eqn:MinPol1}
M_{\varepsilon}(\varepsilon)=\varepsilon^3+c_2\varepsilon^2+c_1\varepsilon-1=0
\end{equation}
with \(c_1,c_2\in\mathbb{Z}\) and we must compute some powers of \(\varepsilon\).
We have \(\varepsilon=1-3r\theta\),
\(\varepsilon^2=(1-3r\theta)^2=1-6r\theta+9r^2\theta^2\), and
\(\varepsilon^3=\varepsilon\cdot\varepsilon^2=
(1-3r\theta)\cdot(1-6r\theta+9r^2\theta)=(1-27r^3b)+(81r^4b-9r)\theta+27r^2\theta^2\).
Substitution into Equation
\eqref{eqn:MinPol1}
yields
\begin{equation}
\label{eqn:MinPol2}
\begin{aligned}
M_{\varepsilon}(\varepsilon) &= (1-27r^3b)+(81r^4b-9r)\theta+27r^2\theta^2 + c_2-6rc_2\theta+9r^2c_2\theta^2 + c_1-3rc_1\theta - 1 \\
&= (1-27r^3b+c_2+c_1-1)+(81r^4b-9r-6rc_2-3rc_1)\theta+(27r^2+9r^2c_2)\theta^2=0.
\end{aligned}
\end{equation}
Consequently, we first obtain \(9r^2(c_2+3)=0\) and thus \(c_2=-3\),
second \(81r^4b-9r+18r-3rc_1=3r(-c_1+3+27r^3b)=0\) and thus \(c_1=3+27r^3b\),
which is confirmed by the third condition \(-27r^3b-3+c_1=0\).
Eventually, we have the \textbf{minimal polynomials}
\begin{equation}
\label{eqn:MinPol3}
M_{\varepsilon}(X)=X^3-3X^2+3(9r^3b+1)X-1 \quad \text{ and } \quad M_{\varepsilon^{-1}}(X)=X^3-3(9r^3b+1)X^2+3X-1.
\end{equation}

Now we apply the results of Nagell
\cite[Satz 22, p. 59]{Ng1930}
and Louboutin
\cite[Theorem 3, pp. 193--194]{Lb2008}
to the minimal polynomials
of our units \(0<\varepsilon<1\) and \(\varepsilon^{-1}>1\) for \(b>1\)
and of the units \(0<\theta<1\) and \(\theta^{-1}>1\) for \(b=1\).

For \(b>1\),
\(\varepsilon\) and \(\varepsilon^{-1}\) are \textbf{fundamental units
of the suborder} \(\mathbb{Z}\lbrack\varepsilon\rbrack\),
since in the nine exceptional sporadic cases
never both, the quadratic and linear coefficient are divisible by \(3\),
and in the exceptional infinite family
it can never happen that \(2r=-3\), resp. \(-2M=3\).
However, since \(\varepsilon=1-3r\theta\) and thus \(\theta=(1-\varepsilon)/3r\),
the index of the suborder in the maximal order is
\(i=(\mathbb{Z}\lbrack\theta\rbrack:\mathbb{Z}\lbrack\varepsilon\rbrack)=27r^3\)
and the discriminant of \(\varepsilon\) is
\(\mathrm{disc}(M_{\varepsilon})=i^2\cdot\mathrm{disc}(M_{\theta})=729r^6\cdot d_P\),
where \(d_P=d_L=-27b^2(4r^3b+1)\).
In Theorem
\ref{thm:VoronoiPeriod},
we shall prove that \(U_L=\langle -1,\varepsilon\rangle\),
i.e., \(\varepsilon\) is a FU of \(L\).

For \(b=1\),
\(\theta\) and \(\theta^{-1}\) are \textbf{fundamental units
of the maximal order} \(\mathcal{O}_L=\mathbb{Z}\lbrack\theta\rbrack\),
because the minimal polynomials
\(M_{\theta}(X)=P(X)=X^3+3rX-1\) and \(M_{\theta^{-1}}(X)=X^3-3rX^2-1\)
do not occur among the exceptional polynomials of Nagell and Louboutin,
for which never a quadratic or linear coefficient vanishes.
\end{proof}


\begin{lemma}
\label{lem:Minima}
If \(\mathcal{O}\subset\mathcal{O}_L\) is a suborder
of the maximal order of a simply real cubic number field \(L\)
and \(\alpha\in\mathcal{O}\) is an algebraic integer
satisfying the condition
\(\lvert\mathrm{N}_{L/\mathbb{Q}}(\alpha)\rvert
<\sqrt[4]{\lvert\mathrm{disc}(\mathcal{O})\rvert/27}\),
then \(\alpha\in\mathrm{Min}(\mathcal{O})\) is a \textbf{lattice minimum}
of the suborder \(\mathcal{O}\),
that is, \((\exists\,i\in\mathbb{Z})\,\alpha=\nu^{i}(1)\).
\end{lemma}

\begin{proof}
The basic idea for this estimate goes back to A. Brentjes.
It was successively improved by G. W. Dueck and H. C. Williams.
See H. C. Williams
\cite[p. 649]{Wi1985}.
\end{proof}


\noindent
The following theorem provides \textit{infinite} families
of \textit{non-pure simply real cubic} number fields \(L=\mathbb{Q}(\theta)\)
whose \textbf{chain of lattice minima} \((\nu^{i}(1))_{i\in\mathbb{Z}}\) 
in the maximal order \(\mathcal{O}_L\) contains non-unitary
generators of absolute differential principal factorizations
(so-called M1- and M2-\textit{chains}, see below for the definition).
This is a sensation, because such infinite families
were known only for \textit{pure cubic} number fields up to now.

\begin{theorem}
\label{thm:Minimum1}
Let \(\theta\) be the real zero of the monogenic trinomial
\(P(X)=X^3+3rbX-b\in\mathbb{Z}\lbrack X\rbrack\)
with squarefree \(b\in\mathbb{N}\) and \(r\in\mathbb{N}\)
and negative discriminant \(d_P=-27b^2(4r^3b+1)\).
The absolute differential principal factor \(\theta\)
of the simply real cubic number field \(L=\mathbb{Q}(\theta)\)
with discriminant \(d_L=d_P\) and maximal order
\(\mathcal{O}_L=\mathbb{Z}\lbrack\theta\rbrack\)
is a \textbf{lattice minimum} in \(\mathrm{Min}(\mathcal{O}_L)\)
if \(r\ge\sqrt[3]{b/4}\).
\end{theorem}

\begin{proof}
We apply Lemma
\ref{lem:Minima}
to the generator \(\theta\) of an absolute differential principal factor
in the maximal order \(\mathcal{O}_L\) of the non-pure simply real cubic number field
\(L=\mathbb{Q}(\theta)\),
according to item (1) of the Main Theorem
\ref{thm:Main}.
By definition, \(\theta\) is the real zero of a monogenic trinomial
\(P(X)=X^3+3rbX-b\in\mathbb{Z}\lbrack X\rbrack\) with 
squarefree absolute coefficient \(b\in\mathbb{N}\), parameter \(r\in\mathbb{N}\),
and sign \(\sigma=+1\).
Therefore, the equation order is maximal,
\(\mathcal{O}_L=\mathbb{Z}\lbrack\theta\rbrack\),
with discriminant \(\mathrm{disc}(\mathcal{O}_L)=d_L=-27b^2(4r^3b+1)\).
By Theorem
\ref{thm:TrinomialPrinciple},
resp. Corollary
\ref{cor:BSS},
we know that 
\(\mathrm{N}_{L/\mathbb{Q}}(\theta)=b\)
and thus the sufficient (but not necessary) condition of Lemma
\ref{lem:Minima}
is equivalent to
\begin{equation*}
\label{eqn:Minima1}
\begin{aligned}
  b         &< \sqrt[4]{27b^2(4r^3b+1)/27} \\
b^4         &< b^2(4r^3b+1) \\
b^2         &< 4r^3b+1 \\
b^2-4r^3b-1 &< 0.
\end{aligned}
\end{equation*}
The expression on the left hand side of this inequality splits into two factors
\begin{equation*}
\label{eqn:Minima2}
\begin{aligned}
\left(b-\frac{1}{2}\left\lbrack 4r^3-\sqrt{16r^6+4}\right\rbrack\right)\cdot\left(b-\frac{1}{2}\left\lbrack 4r^3+\sqrt{16r^6+4}\right\rbrack\right) &< 0 \\
\left(b-\left\lbrack 2r^3-\sqrt{4r^6+1}\right\rbrack\right)\cdot\left(b-\left\lbrack 2r^3+\sqrt{4r^6+1}\right\rbrack\right) &< 0 \\
\left(b-2r^3\left\lbrack 1-\sqrt{1+\frac{1}{4r^6}}\right\rbrack\right)\cdot\left(b-2r^3\left\lbrack 1+\sqrt{1+\frac{1}{4r^6}}\right\rbrack\right) &< 0.
\end{aligned}
\end{equation*}
Since \(b\) and \(r\) are positive, we must have
\begin{equation*}
\label{eqn:Minima3}
\begin{aligned}
b-2r^3\left\lbrack 1-\sqrt{1+\frac{1}{4r^6}}\right\rbrack>0 \quad &\text{ and } \quad 
b-2r^3\left\lbrack 1+\sqrt{1+\frac{1}{4r^6}}\right\rbrack<0 \\
2r^3\left\lbrack 1-\sqrt{1+\frac{1}{4r^6}}\right\rbrack<b \quad &\text{ and } \quad 
b<2r^3\left\lbrack 1+\sqrt{1+\frac{1}{4r^6}}\right\rbrack,
\end{aligned}
\end{equation*}
where the left inequality is trivial, and the right inequality can be sharpened to \(b\le 4r^3\).
For fixed \(b\), the parameter \(r\) must be bounded from below by \(r\ge\sqrt[3]{b/4}\).
\end{proof}

\noindent
For each fixed value of the absolute coefficient \(b\),
this sufficient (but not necessary) condition admits an \textit{infinitude} of parameters \(r\),
for instance,
\begin{equation*}
\label{eqn:Minima4}
\begin{aligned}
r \ge \sqrt[3]{1} = 1  &\text{ for } b = 4 = 4\cdot 1^3, \\
r \ge \sqrt[3]{8} = 2  &\text{ for } b = 32 = 4\cdot 2^3, \\
r \ge \sqrt[3]{27} = 3 &\text{ for } b = 108 = 4\cdot 3^3, \\
r \ge \sqrt[3]{64} = 4 &\text{ for } b = 256 = 4\cdot 4^3.
\end{aligned}
\end{equation*}


\noindent
We supplement Theorem
\ref{thm:Minimum1}
by a statement about a \textit{second absolute differential principal factor}
which has less probability to be a lattice minimum.

\medskip
\noindent
If \(\ell\) denotes the \textit{primitive period length}
of the \textbf{Voronoi Algorithm}
\cite{Vo1896}
for the system of lattice minima
\(\mathrm{Min}(\mathcal{O}_L)\subset\mathcal{O}_L\)
in the maximal order of the simply real cubic number field
\(L=\mathbb{Q}(\theta)\),
then the two-sided chain
\((\nu^{i}(1))_{i\in\mathbb{Z}}\)
consists precisely of \(\ell\) pairwise disjoint \textit{orbits}
\(\dot{\cup}_{j=1}^{\ell}\,U_L^+\cdot\nu^{j}(1)\)
under the operation of the \textit{positive} unit group \(U_L^+\) on \(\mathrm{Min}(\mathcal{O}_L)\).
One speaks about an M2-chain, resp. M1-chain, resp. M0-chain, 
if there are \(2\), resp. \(1\), resp. \(0\), orbits
of non-unitary differential principal factors in \(\mathrm{Min}(\mathcal{O}_L)\).
(Since \(\nu^{\ell}(1)=\varepsilon\), the trivial \textit{unitary} orbit is always \((\nu^{k\cdot \ell}(1))_{k\in\mathbb{Z}}\).)

\begin{corollary}
\label{cor:Minimum2}
Let \(\theta\) be the real zero of the \textbf{monogenic} trinomial
\(P(X)=X^3+3rbX-b\in\mathbb{Z}\lbrack X\rbrack\)
with parameter \(r\in\mathbb{N}\) and square free coefficient \(b\ge 2\).
The two-sided chain
\((\nu^i(1))_{i\in\mathbb{Z}}\)
of lattice minima \(\nu^i(1)\in\mathrm{Min}(\mathcal{O}_L)\)
in the discrete geometric Minkowski-image of the maximal order \(\mathcal{O}_L\)
(a complete lattice in Euclidean space \(\mathbb{R}^3\))
of the simply real cubic number field \(L=\mathbb{Q}(\theta)\) is
\begin{itemize}
\item
a guaranteed \(\mathrm{M}2\)-chain, if 
\(r\ge b\cdot\sqrt[3]{b^2/4}\left(>\sqrt[3]{b/4}\right)\), and
\item
at least an \(\mathrm{M}1\)-chain (possibly even an \(\mathrm{M}2\)-chain), if
\(r\ge\sqrt[3]{b/4}\).
\end{itemize}
\end{corollary}

\begin{proof}
Similarly as in the proof of Theorem
\ref{thm:Minimum1},
we apply Lemma
\ref{lem:Minima}
to the second generator \(\theta^2\) of an absolute differential principal factor
in the maximal order \(\mathcal{O}_L\).
By Theorem
\ref{thm:TrinomialPrinciple},
resp. Corollary
\ref{cor:BSS},
we have
\(\mathrm{N}_{L/\mathbb{Q}}(\theta^2)=b^2\)
and thus the condition of Lemma
\ref{lem:Minima}
is equivalent to
\begin{equation*}
\label{eqn:Minimum2}
\begin{aligned}
  b^2       &< \sqrt[4]{27b^2(4r^3b+1)/27} \\
b^8         &< b^2(4r^3b+1) \\
b^6         &< 4r^3b+1.
\end{aligned}
\end{equation*}
Now we abstain from an exact solution of this inequality
and sharpen immediately to
\(b^6\le 4r^3b\), resp. \(b^5\le 4r^3\).
Consequently,
for fixed \(b\), the parameter \(r\) must be bounded from below by
\(r\ge\sqrt[3]{b^5/4}=b\cdot\sqrt[3]{b^2/4}=b\cdot\sqrt[3]{(b/2)^2}\).
\end{proof}

\noindent
Again, for each fixed value of the absolute coefficient \(b\),
this sufficient (but not necessary) condition admits an \textit{infinitude} of parameters \(r\)
for warranted \(\mathrm{M}2\)-chains,
for instance,
\begin{equation*}
\label{eqn:Minimum2Exm}
\begin{aligned}
r \ge 2\cdot\sqrt[3]{1} = 2\cdot 1 = 2            &\text{ for } b = 2 = 2\cdot 1^3, \\
r \ge 16\cdot\sqrt[3]{64} = 16\cdot 4 = 64        &\text{ for } b = 16 = 2\cdot 2^3, \\
r \ge 54\cdot\sqrt[3]{729} = 54\cdot 9 = 486      &\text{ for } b = 54 = 2\cdot 3^3, \\
r \ge 128\cdot\sqrt[3]{4096} = 128\cdot 16 = 2048 &\text{ for } b = 128 = 2\cdot 4^3.
\end{aligned}
\end{equation*}


\noindent
Since the criterion in Lemma
\ref{lem:Minima}
is a coarse sufficient but not necessary condition,
we try a direct approach to a sharp criterion
for the containment of the absolute DPF \(\theta\)
in \(\mathrm{Min}(\mathcal{O}_L)\):
If \(\theta\not\in\mathrm{Min}(\mathcal{O}_L)\) is not a lattice minimum in the maximal order
\(\mathcal{O}_L=\mathbb{Z}\lbrack\theta\rbrack\) of the \textit{simply real} cubic number field \(L=\mathbb{Q}(\theta)\),
then there exists an algebraic integer \(\nu\in\mathcal{O}_L\) satisfying the inequalities
\(0<\nu<\theta<1\) and \(1<\nu^\prime\nu^{\prime\prime}<\theta^\prime\theta^{\prime\prime}\),
i.e., the norm cylinder of \(\theta\) contains a non trivial lattice point \(\nu\).
Now the decisive clue is the investigation of the \textbf{crucial element}
\(\phi:=\mathrm{N}(\theta)\cdot\nu/\theta\in L\).

\begin{proposition}
\label{prp:Crucial}
\(\phi\) is an algebraic integer in \(\mathcal{O}_L=\mathbb{Z}\lbrack\theta\rbrack\),
whose norm \(\mathrm{N}(\phi)=b^2\cdot\mathrm{N}(\nu)\) is divisible by \(b^2\).
Due to the monogeneity, \(\phi=x+y\theta+z\theta^2\) with integers \(x,y,z\in\mathbb{Z}\),
such that \(b\mid x\), \(b\mid y\), and \(z^3\equiv\mathrm{N}(\nu)\,(\mathrm{mod}\,b)\).
The \textbf{complement} and the \textbf{norm} of \(\phi\) are given explicitly by
\begin{equation}
\label{eqn:Crucial}
\begin{aligned}
\phi^\prime\phi^{\prime\prime} &= (x^2+3rby^2+9r^2b^2z^2-byz-6rbxz)+(bz^2-xy)\theta+(y^2+3rbz^2-xz)\theta^2, \\
\mathrm{N}(\phi)               &= x^3+by^3+b^2z^3+3rbxy^2+3rb^2yz^2+9r^2b^2xz^2-6rbx^2z-3bxyz.
\end{aligned}
\end{equation}
The bounds \(0<\phi<b\) and \(\phi^\prime\phi^{\prime\prime}<b^2\), i.e. \(\lvert\phi^\prime\rvert=\lvert\phi^{\prime\prime}\rvert<b\), imply three estimates,
\begin{equation}
\label{eqn:Estimates}
\lvert x-2rbz\rvert<b, \quad \lvert 3z-6ry\rvert<1+2\sqrt{3rb+1}, \quad \lvert 6r^2bz+y-2rx\rvert<1+2rb.
\end{equation}
\end{proposition}

\begin{proof}
Since \(\nu<\theta\) and thus \(\nu/\theta<1\), its \textit{height} is bounded by
\(\phi=\mathrm{N}(\theta)\cdot\nu/\theta<\mathrm{N}(\theta)=b\),
and since \(\nu^\prime\nu^{\prime\prime}<\theta^\prime\theta^{\prime\prime}\)
and thus \(\nu^\prime\nu^{\prime\prime}/(\theta^\prime\theta^{\prime\prime})<1\),
the \textit{square of the radius} is bounded by
\(\phi^\prime\phi^{\prime\prime}=\mathrm{N}(\theta)^2\cdot\nu^\prime\nu^{\prime\prime}/(\theta^\prime\theta^{\prime\prime})
<\mathrm{N}(\theta)^2=b^2\).
The element
\(\phi=\theta\theta^\prime\theta^{\prime\prime}\cdot\nu/\theta=\theta^\prime\theta^{\prime\prime}\cdot\nu\)
is an algebraic \textit{integer} in \(\mathcal{O}_L\),
since the quadratic elementary symmetric polynomial
\(\theta\theta^{\prime}+\theta^\prime\theta^{\prime\prime}+\theta^{\prime\prime}\theta\)
is the linear coefficient \(3rb\) of the minimal polynomial \(P(X)=X^3+3rbX-b\) of \(\theta\),
and the trace \(\theta+\theta^\prime+\theta^{\prime\prime}=0\) vanishes, whence
\(\theta^\prime\theta^{\prime\prime}=3rb-\theta\theta^{\prime}-\theta^{\prime\prime}\theta
=3rb-\theta\cdot (\theta^{\prime}+\theta^{\prime\prime})=3rb-\theta\cdot (-\theta)
=3rb+\theta^2\in\mathcal{O}_L\).
The \textit{norm} of \(\phi\), \(\mathrm{N}(\phi)=\mathrm{N}(\theta)^3\cdot\mathrm{N}(\nu)/\mathrm{N}(\theta)
=\mathrm{N}(\theta)^2\cdot\mathrm{N}(\nu)=b^2\cdot\mathrm{N}(\nu)\), is divisible by \(b^2\).
Since \(\mathcal{O}_L=\mathbb{Z}\lbrack\theta\rbrack\) is \textit{monogenic},
we have \(\phi=x+y\theta+z\theta^2\) with integers \(x,y,z\in\mathbb{Z}\).
The \textit{complement} (i.e., the square of the radius) of \(\phi\) is given by
\(\phi^\prime\phi^{\prime\prime}=(x+y\theta^\prime+z(\theta^\prime)^2)\cdot(x+y\theta^{\prime\prime}+z(\theta^{\prime\prime})^2)\)
\begin{equation*}
\label{eqn:complement}
\begin{aligned}
&= x^2+y^2\theta^\prime\theta^{\prime\prime}+z^2(\theta^\prime\theta^{\prime\prime})^2
+xy(\theta^{\prime}+\theta^{\prime\prime})
+yz(\theta^\prime(\theta^{\prime\prime})^2+\theta^{\prime\prime}(\theta^\prime)^2)
+xz((\theta^\prime)^2+(\theta^{\prime\prime})^2) \\
&= x^2+y^2(3rb+\theta^2)+z^2(9r^2b^2+6rb\theta^2+b\theta-3rb\theta^2)
-xy\theta-byz-xz(\theta^2+6rb) \\
&= (x^2+3rby^2+9r^2b^2z^2-byz-6rbxz)+(bz^2-xy)\theta+(y^2+3rbz^2-xz)\theta^2,
\end{aligned}
\end{equation*}
since \(\theta^\prime(\theta^{\prime\prime})^2+\theta^{\prime\prime}(\theta^\prime)^2
=\theta^\prime\theta^{\prime\prime}(\theta^{\prime\prime}+\theta^\prime)
=(3rb+\theta^2)\cdot(-\theta)=-3rb\theta-\theta^3=-3rb\theta-(b-3rb\theta)=-b\)
and \((\theta^\prime)^2+(\theta^{\prime\prime})^2
=(\theta^\prime+\theta^{\prime\prime})^2-2\theta^\prime\theta^{\prime\prime}
=(-\theta)^2-2(3rb+\theta^2)=\theta^2-6rb-2\theta^2=-(\theta^2+6rb)\).
Now we can determine the \textit{norm} \(\mathrm{N}(\phi)\in\mathbb{Z}\) of \(\phi\) explicitly,
neglecting all terms with \(\theta\) and \(\theta^2\):
\(\mathrm{N}(\phi)=\phi\cdot\phi^\prime\phi^{\prime\prime}
=x(x^2+3rby^2+9r^2b^2z^2-byz-6rbxz)+y(y^2+3rbz^2-xz)b+z(bz^2-xy)b\)
\(=x^3+3rbxy^2+9r^2b^2xz^2-bxyz-6rbx^2z+by^3+3rb^2yz^2-bxyz+b^2z^3-bxyz\) \\
\(=x^3+by^3+b^2z^3+3rbxy^2+3rb^2yz^2+9r^2b^2xz^2-6rbx^2z-3bxyz\). \\
The divisibility 
\(\mathrm{N}(\phi)=b^2\cdot\mathrm{N}(\nu)\)
by \(b^2\) implies \(b\mid x^3\) and thus, since \(b\) is squarefree, \(x=b\cdot\tilde{x}\) and
\(b\cdot\mathrm{N}(\nu)=b^2\tilde{x}^3+y^3+bz^3+3rb\tilde{x}y^2+3rbyz^2+9r^2b^2\tilde{x}z^2-6rb^2\tilde{x}^2z-3b\tilde{x}yz\).
Thus \(y=b\cdot\tilde{y}\) and
\(\mathrm{N}(\nu)=b\tilde{x}^3+b^2\tilde{y}^3+z^3+3rb^2\tilde{x}\tilde{y}^2+3rb\tilde{y}z^2+9r^2b\tilde{x}z^2-6rb\tilde{x}^2z-3b\tilde{x}\tilde{y}z\),
whence \(\mathrm{N}(\nu)\equiv z^3\,(\mathrm{mod}\,b)\).
Finally, we derive three estimates, using the triangle inequality: Firstly, \\
\(\phi+\phi^\prime+\phi^{\prime\prime}
=3x+y(\theta+\theta^\prime+\theta^{\prime\prime})+z(\theta^2+(\theta^\prime)^2+(\theta^{\prime\prime})^2)
=3x+y\cdot 0+z\cdot(\theta^2-6rb-\theta^2)=3x-6rb\) and thus
\(\lvert 3x-6rbz\rvert=\lvert\phi+\phi^\prime+\phi^{\prime\prime}\rvert
\le\lvert\phi\rvert+\lvert\phi^\prime\rvert+\lvert\phi^{\prime\prime}\rvert<3b\),
resp. \(\lvert x-2rbz\rvert<b\). Secondly, \\
\(\theta\phi+\theta^\prime\phi^\prime+\theta^{\prime\prime}\phi^{\prime\prime}
=x(\theta+\theta^\prime+\theta^{\prime\prime})+y(\theta^2+(\theta^\prime)^2+(\theta^{\prime\prime})^2)
+z(\theta^3+(\theta^\prime)^3+(\theta^{\prime\prime})^3)=3bz-6rby\), \\
since \(\theta+\theta^\prime+\theta^{\prime\prime}=0\),
\(\theta^2+(\theta^\prime)^2+(\theta^{\prime\prime})^2=-6rb\),
\(\theta^3=b-3rb\theta\), and
\((\theta^\prime)^3+(\theta^{\prime\prime})^3
=(\theta^\prime+\theta^{\prime\prime})^3-3(\theta^\prime(\theta^{\prime\prime})^2+\theta^{\prime\prime}(\theta^\prime)^2)
=(-\theta)^3-3(-b)=3rb\theta-b+3b=3rb\theta+2b\), whence
\(\theta^3+(\theta^\prime)^3+(\theta^{\prime\prime})^3=b-3rb\theta+3rb\theta+2b=3b\). Thus, 
\(\lvert 3bz-6rby\rvert=\lvert\theta\phi+\theta^\prime\phi^\prime+\theta^{\prime\prime}\phi^{\prime\prime}\rvert
\le\lvert\theta\rvert\lvert\phi\rvert
+\lvert\theta^\prime\rvert\lvert\phi^\prime\rvert
+\lvert\theta^{\prime\prime}\rvert\lvert\phi^{\prime\prime}\rvert
<1\cdot b+\sqrt{3rb+1}\cdot b+\sqrt{3rb+1}\cdot b=b(1+2\sqrt{3rb+1})\),
resp. \(\lvert 3z-6ry\rvert<1+2\sqrt{3rb+1}\),
since \(\lvert\theta\rvert<1\) and
\(\lvert\theta^\prime\rvert^2=\lvert\theta^{\prime\prime}\rvert^2
=\lvert\theta^\prime\theta^{\prime\prime}\rvert=\lvert 3rb+\theta^2\rvert
\le 3rb+\lvert\theta\rvert^2<3rb+1\).
Eventually, \\
\(\theta^2\phi+(\theta^\prime)^2\phi^\prime+(\theta^{\prime\prime})^2\phi^{\prime\prime}
=x(\theta^2+(\theta^\prime)^2+(\theta^{\prime\prime})^2)+y(\theta^3+(\theta^\prime)^3+(\theta^{\prime\prime})^3)
+z(\theta^4+(\theta^\prime)^4+(\theta^{\prime\prime})^4)=18r^2b^2z+3by-6rbx\),
since \(\theta^4=b\theta-3rb\theta^2\)
and \((\theta^\prime)^4+(\theta^{\prime\prime})^4
=(\theta^\prime+\theta^{\prime\prime})^4
-4(\theta^\prime)^3\theta^{\prime\prime}-6(\theta^\prime)^2(\theta^{\prime\prime})^2-4\theta^\prime(\theta^{\prime\prime})^3
=(-\theta)^4-4\theta^\prime\theta^{\prime\prime}((\theta^\prime)^2+(\theta^{\prime\prime})^2)-6(9r^2b^2+6rb\theta^2+\theta^4)
=b\theta-3rb\theta^2-4(3rb+\theta^2)\cdot\) \((-6rb-\theta^2)-6(9r^2b^2+6rb\theta^2+b\theta-3rb\theta^2)
=b\theta-3rb\theta^2+4(18r^2b^2+9rb\theta^2+b\theta-3rb\theta^2)-54r^2b^2-6b\theta-18rb\theta^2
=18r^2b^2-b\theta+3rb\theta^2\). Thus
\(\lvert 18r^2b^2z+3by-6rbx\rvert=\lvert\theta^2\phi+(\theta^\prime)^2\phi^\prime+(\theta^{\prime\prime})^2\phi^{\prime\prime}\rvert
\le\lvert\theta\rvert^2\lvert\phi\rvert+\lvert\theta^\prime\rvert^2\lvert\phi^\prime\rvert+\lvert\theta^{\prime\prime}\rvert^2\lvert\phi^{\prime\prime}\rvert
<1\cdot b+(3rb+1)\cdot b+(3rb+1)\cdot b=b(1+2+6rb)=3b(1+2rb)\),
resp. \(\lvert 6r^2bz+y-2rx\rvert<1+2rb\).
\end{proof}


\noindent
We combine the inequalities in Formula
\eqref{eqn:Estimates}
with divisibility properties,
which have not been exploited yet.

\begin{lemma}
\label{lem:Crucial}
Let \(b\ge 2\) be a squarefree integer and \(r\ge 1\) be a positive integer.
A triplet of integers \((x,y,z)\in\mathbb{Z}^3\)
which satisfies the estimates
\begin{equation}
\label{eqn:EstimatesClone}
\lvert x-2rbz\rvert<b, \quad \lvert 3z-6ry\rvert<1+2\sqrt{3rb+1}, \quad \lvert 6r^2bz+y-2rx\rvert<1+2rb,
\end{equation}
under the divisibility constraints \(b\mid x\) and \(b\mid y\),
must be trivial, that is
\begin{equation}
\label{eqn:Trivial}
(x,y,z)=(0,0,0).
\end{equation}
\end{lemma}

\begin{proof}
Since \(b\mid x\), there exists an integer \(\tilde{x}\) such that \(x=b\cdot\tilde{x}\).
The \textit{first} estimate
\(\lvert b\tilde{x}-2rbz\rvert<b\), divided by \(b\), yields
\(\lvert \tilde{x}-2rz\rvert<1\) and thus \(\tilde{x}=2rz\).
Since \(b\mid y\), there exists an integer \(\tilde{y}\) such that \(y=b\cdot\tilde{y}\).
Substitution into the \textit{third} estimate shows
\(\lvert 6r^2bz+b\tilde{y}-2r\cdot 2rbz\rvert<1+2rb\), resp.
\(\lvert 2r^2z+\tilde{y}\rvert=\lvert 6r^2z+\tilde{y}-4r^2z\rvert<2r+\frac{1}{b}\).
Since we seek integer solutions, this can be sharpened to
\(-2r\le 2r^2z+\tilde{y}\le 2r\), resp.
\(-2r-2r^2z\le\tilde{y}\le 2r-2r^2z\) and thus
\(-2r(1+rz)\le\tilde{y}\le 2r(1-rz)\).
This enforces \(1+rz>0\), \(1-rz>0\),
and consequently \(-1<rz<1\),
which can only be satisfied by \(rz=0\), and thus \(z=0\).
There remains a possible range \(-2r\le\tilde{y}\le 2r\),
while the other two components are already narrowed down to
\(\tilde{x}=z=0\).
Finally we consider the \textit{second} estimate, which degenerates to
\(6rb\lvert\tilde{y}\rvert<1+2\sqrt{3rb+1}\), resp.
\(\lvert\tilde{y}\rvert<\frac{1}{6rb}+\frac{\sqrt{3rb+1}}{3rb}\).
The sum of the two fractions on the right hand side, however,
is very small, since \(\frac{1}{6rb}\le\frac{1}{12}\) and
\(\frac{\sqrt{3rb+1}}{3rb}<\frac{\sqrt{3rb+3rb}}{3rb}
=\frac{\sqrt{2}}{\sqrt{3rb}}\le\frac{\sqrt{2}}{\sqrt{6}}=\frac{1}{\sqrt{3}}<\frac{7}{12}\).
Eventually,
the solution of the resulting inequality \(\lvert\tilde{y}\rvert<\frac{8}{12}=\frac{2}{3}\)
is trivial, \(\tilde{y}=0\).
\end{proof}


\noindent
This auxiliary result has the following striking consequence,
which immediately sharpens Theorem
\ref{thm:Minimum1}
and the second item of Corollary
\ref{cor:Minimum2},
without the requirement of any further proof.

\begin{theorem}
\label{thm:Unconditional}
Let \(\theta\) be the real zero of the \textbf{monogenic} trinomial
\(P(X)=X^3+3rbX-b\in\mathbb{Z}\lbrack X\rbrack\)
with parameter \(r\in\mathbb{N}\), square free coefficient \(b\ge 2\),
and negative discriminant \(d_P=-27b^2(4r^3b+1)<0\).

Then the non-unitary absolute differential principal factor \(\theta\)
of the simply real cubic number field \(L=\mathbb{Q}(\theta)\)
with discriminant \(d_L=d_P\) and maximal order
\(\mathcal{O}_L=\mathbb{Z}\lbrack\theta\rbrack\)
is \textbf{always a lattice minimum} in \(\mathrm{Min}(\mathcal{O}_L)\),
and thus can be found with the Voronoi algorithm.

The two-sided chain
\((\nu^i(1))_{i\in\mathbb{Z}}\)
of lattice minima \(\nu^i(1)\in\mathrm{Min}(\mathcal{O}_L)\)
in the discrete geometric Minkowski-image of \(\mathcal{O}_L\)
is \textbf{always at least an \(\mathrm{M}1\)-chain} (possibly even an \(\mathrm{M}2\)-chain).
\end{theorem}


\begin{remark}
\label{rmk:Crucial}
The proof of Proposition
\ref{prp:Crucial}
was conducted under the assumption that the real zero \(\theta\)
of the polynomial \(X^3+3rbX-b\in\mathbb{Z}\lbrack X\rbrack\)
(with \(\sigma=+1\))
is contained in the open interval \(0<\theta<1\).
Actually,
this is true more generally for any real zero \(\theta\)
of a trinomial \(P(X)=X^p+vbX-b\in\mathbb{Z}\lbrack X\rbrack\)
with positive integers \(b\), \(v\),
and any \textit{odd} prime number \(p\).
If we had \(\theta\le 0\), then \(\theta^p\le 0\) and thus
\(\theta^p+vb\theta-b\le -b<0\), in contradiction to \(P(\theta)=0\).
On the other hand, if \(\theta\ge 1\), then
\(\theta^p+vb\theta-b\ge 1+vb-b\ge 1+b-b=1>0\), which is also a contradiction.
\end{remark}


\noindent
We proceed with an \textit{explicit parametrization}
of the \textit{first primitive period} of the Voronoi chain
in a \textit{simply real} cubic field,
which confirms Theorem
\ref{thm:Unconditional}
and improves Corollary
\ref{cor:Minimum2}.

\begin{theorem}
\label{thm:VoronoiPeriod}
\textbf{(Parametrized Period.)}
The period length \(\ell\) of the Voronoi algorithm for
the maximal order \(\mathcal{O}_L\) of a
simply real cubic number field \(L=\mathbb{Q}(\theta)\),
generated by a \textbf{monogenic} trinomial \(P(X)=X^3+3rbX-b\)
with parameter \(r\ge 1\) and squarefree absolute coefficient \(b\ge 1\),
can take only two values,
either \(\ell=2\) when \(r<L(b)\) and the chain \(\mathrm{Min}(\mathcal{O}_L)\) is an \(\mathrm{M}1\)-chain,
\begin{equation}
\label{eqn:PL2}
\begin{aligned}
\nu^0(1)=1,      & \qquad & \overbrace{\nu^1(1)=\theta}^{\mathrm{DPF}}, & \qquad & \overbrace{\nu^2(1)=1-3r\theta=\varepsilon_0}^{\mathrm{FU}}, \\
\mathrm{N}(1)=1, & \qquad & \mathrm{N}(\theta)=b,              & \qquad & \mathrm{N}(\varepsilon_0)=1,
\end{aligned}
\end{equation}
or \(\ell=3\) when \(r\ge L(b)\) and the Voronoi chain \(\mathrm{Min}(\mathcal{O}_L)\) is an \(\mathrm{M}2\)-chain,
\begin{equation}
\label{eqn:PL3}
\begin{aligned}
\nu^0(1)=1,      & \qquad & \overbrace{\nu^1(1)=\theta}^{1\mathrm{st\ DPF}}, & \qquad & \overbrace{\nu^2(1)=\theta^2}^{2\mathrm{nd\ DPF}}, & \qquad & \overbrace{\nu^3(1)=1-3r\theta=\varepsilon_0}^{\mathrm{FU}}, \\
\mathrm{N}(1)=1, & \qquad & \mathrm{N}(\theta)=b,              & \qquad & \mathrm{N}(\theta^2)=b^2,            & \qquad & \mathrm{N}(\varepsilon_0)=1,
\end{aligned}
\end{equation}
in dependence on a \textbf{sharp lower bound} \(L(b)\) for the parameter \(r\),
for each fixed coefficient \(b\ge 2\).
For \(b=1\), however, the period length takes the minimal value \(\ell=1\) and the \(\mathrm{M}0\)-chain is
\begin{equation}
\label{eqn:PL1}
\begin{aligned}
\nu^0(1)=1,      & \qquad & \overbrace{\nu^1(1)=\theta=\varepsilon_0}^{\mathrm{FU}}, & \qquad \Biggl\lbrack & \nu^2(1)=\theta^2=\varepsilon_0^2, & \qquad & \nu^3(1)=1-3r\theta=\varepsilon_0^3,\ \Biggr\rbrack \\
\mathrm{N}(1)=1, & \qquad & \mathrm{N}(\varepsilon_0)=1,                             & \qquad \Biggl\lbrack & \mathrm{N}(\varepsilon_0^2)=1,     & \qquad & \mathrm{N}(\varepsilon_0^3)=1.\ \Biggr\rbrack
\end{aligned}
\end{equation}
In the cases with \(\ell\in\lbrace 2,3\rbrace\), all non-unitary lattice minima are \textbf{differential principal factors}.
\end{theorem}

\begin{proof}
We use the Voronoi algorithm as described in Williams/Zarnke
\cite[pp. 442--446]{WiZa1972}.
On principle, we could also refer to Delone/Faddeev
\cite[pp. 282--286]{DeFa1964},
but these 5 pages contain 14 typographical errors.
See also
\cite{BWZ1971,HaWi2018}.
First, we consider the \textit{regular} situation
\(b\ge 2\) with usual period length \(\ell=3\).
Formally, the scenario \(b=1\) with \(\ell=1\)
is contained as a special case in Formula
\eqref{eqn:PL3},
indicated by brackets in Formula
\eqref{eqn:PL1}. \\
\textbf{Initialization}: We start with
\begin{equation}
\label{eqn:Initialization}
(1,\theta,\theta^2)\cdot
\begin{pmatrix}
1 & 0 & 0 \\
0 & 1 & 0 \\
0 & 0 & 1
\end{pmatrix}
=(1,\theta,\theta^2), \text{ the \textbf{monogenic} integral basis of } \mathcal{O}_L=\mathbb{Z}\lbrack\theta\rbrack.
\end{equation}
\(0\)th \textbf{Reduction}: We execute the steps from 1) to 6) as described in
\cite[pp. 444--446]{WiZa1972},
and we get
\begin{equation}
\label{eqn:Reduction0}
(1,\theta,\theta^2)\cdot
\begin{pmatrix}
1 & 0 & 1 \\
0 & 1 & -1 \\
0 & 0 & 1
\end{pmatrix}
=(1,\vartheta^0_g,\vartheta^0_h), \text{ the \(0\)th reduced Voronoi basis.}
\end{equation}
Thus, the first non-trivial lattice minimum is
\(\nu^1(1)=\vartheta^0_g=\theta\)
and its supplement is
\(\vartheta^0_h=1-\theta+\theta^2\). \\
\(1\)st \textbf{Division}:
Instead of dividing by \(\vartheta^0_g\),
we multiply by the quotient \textbf{complement by norm}, since
\(1/\vartheta^0_g
=(\vartheta^0_g)^\prime(\vartheta^0_g)^{\prime\prime}/\vartheta^0_g(\vartheta^0_g)^\prime(\vartheta^0_g)^{\prime\prime}
=(\vartheta^0_g)^\prime(\vartheta^0_g)^{\prime\prime}/\mathrm{N}(\vartheta^0_g)\).
According to Formula
\eqref{eqn:Crucial},
\((\vartheta^0_g)^\prime(\vartheta^0_g)^{\prime\prime}=\theta^\prime\theta^{\prime\prime}=3rb+\theta^2\)
and \(\mathrm{N}(\vartheta^0_g)=\mathrm{N}(\theta)=b\).
Thus, \(1/\vartheta^0_g=1/\theta=(3rb+\theta^2)/b\) and
\(\vartheta^0_h/\vartheta^0_g=(1-\theta+\theta^2)(3rb+\theta^2)/b
=(3rb-b+b\theta+\theta^2)/b\),
and we get the \(1\)st divided basis
\begin{equation}
\label{eqn:Division1}
(1,\theta,\theta^2)\cdot\frac{1}{b}
\begin{pmatrix}
b & 3rb-b & 3rb \\
0 & b     & 0 \\
0 & 1     & 1
\end{pmatrix}
=(1,\vartheta^0_h/\vartheta^0_g,1/\vartheta^0_g).
\end{equation}
\(1\)st \textbf{Reduction}: Again, we execute the steps from 1) to 6),
and we get
\begin{equation}
\label{eqn:Reduction1}
(1,\theta,\theta^2)\cdot\frac{1}{b}
\begin{pmatrix}
b & 0 & 0 \\
0 & b & 0 \\
0 & 0 & 1
\end{pmatrix}
=(1,\vartheta^1_g,\vartheta^1_h), \text{ the \(1\)st reduced Voronoi basis.}
\end{equation}
Thus, the second \textbf{relative} lattice minimum is
\(\vartheta^1_g=b\theta/b=\theta\)
and its supplement is
\(\vartheta^1_h=\theta^2/b\).
The second \textbf{absolute} lattice minimum is the product
\(\nu^2(1)=\vartheta^0_g\vartheta^1_g=\theta\cdot\theta=\theta^2\). \\
\(2\)nd \textbf{Division}:
According to Formula
\eqref{eqn:Crucial},
\((\vartheta^1_g)^\prime(\vartheta^1_g)^{\prime\prime}=\theta^\prime\theta^{\prime\prime}=3rb+\theta^2\)
and \(\mathrm{N}(\vartheta^1_g)=\mathrm{N}(\theta)=b\).
Thus, \(1/\vartheta^1_g=1/\theta=(3rb+\theta^2)/b\) and
\(\vartheta^1_h/\vartheta^1_g=\theta^2(3rb+\theta^2)/b^2
=b\theta/b^2=\theta/b\),
and we get the \(2\)nd divided basis
\begin{equation}
\label{eqn:Division2}
(1,\theta,\theta^2)\cdot\frac{1}{b}
\begin{pmatrix}
b & 0 & 3rb \\
0 & 1 & 0 \\
0 & 0 & 1
\end{pmatrix}
=(1,\vartheta^1_h/\vartheta^1_g,1/\vartheta^1_g).
\end{equation}
\(2\)nd \textbf{Reduction}: Again, we execute the steps from 1) to 6),
and we get
\begin{equation}
\label{eqn:Reduction2}
(1,\theta,\theta^2)\cdot\frac{1}{b}
\begin{pmatrix}
b & 0 & b \\
0 & 1 & -1 \\
0 & 0 & 1
\end{pmatrix}
=(1,\vartheta^2_g,\vartheta^2_h), \text{ the \(2\)nd reduced Voronoi basis.}
\end{equation}
Thus, the third \textbf{relative} lattice minimum is
\(\vartheta^2_g=\theta/b\)
and its supplement is
\(\vartheta^2_h=(b-\theta+\theta^2)/b\).
The third \textbf{absolute} lattice minimum is the product
\(\nu^3(1)=\vartheta^0_g\vartheta^1_g\vartheta^2_g=\theta\cdot\theta\cdot\theta/b
=\theta^3/b=1-3r\theta=\varepsilon_0\),
that is the desired \textbf{fundamental unit} (FU).
Although the primitive period of lattice minima with length \(\ell=3\) is found,
yet another division and reduction step is required in order
to finish the algorithm with a tidy break off condition. \\
\(3\)rd \textbf{Division}:
According to Formula
\eqref{eqn:Crucial},
\((\vartheta^2_g)^\prime(\vartheta^2_g)^{\prime\prime}=\theta^\prime\theta^{\prime\prime}/b^2=(3rb+\theta^2)/b^2\)
and \(\mathrm{N}(\vartheta^2_g)=\mathrm{N}(\theta/b)=b/b^3=1/b^2\).
Thus, \(1/\vartheta^2_g=b/\theta=((3rb+\theta^2)/b^2)/(1/b^2)=3rb+\theta^2\) and
\(\vartheta^2_h/\vartheta^2_g=(b-\theta+\theta^2)(3rb+\theta^2)/b
=(3rb^2-b+b\theta+b\theta^2)/b=3rb-1+\theta+\theta^2\),
and the \(3\)rd divided basis is
\begin{equation}
\label{eqn:Division3}
(1,\theta,\theta^2)\cdot
\begin{pmatrix}
1 & 3rb-1 & 3rb \\
0 & 1 & 0 \\
0 & 1 & 1
\end{pmatrix}
=(1,\vartheta^2_h/\vartheta^2_g,1/\vartheta^2_g).
\end{equation}
\(3\)rd \textbf{Reduction}: For a last time, we execute the steps from 1) to 6),
and we get
\begin{equation}
\label{eqn:Reduction3}
(1,\theta,\theta^2)\cdot
\begin{pmatrix}
1 & 0 & 1 \\
0 & 1 & -1 \\
0 & 0 & 1
\end{pmatrix}
=(1,\vartheta^3_g,\vartheta^3_h), \text{ the \(3\)rd reduced Voronoi basis.}
\end{equation}
This is the same lattice as in Formula
\eqref{eqn:Reduction0}
and the Voronoi algorithm can be terminated.
We point out that the main denominator in Formulae
\eqref{eqn:Division1}--\eqref{eqn:Reduction2}
is \(b\), whereas it is \(1\) in Formulae
\eqref{eqn:Initialization},
\eqref{eqn:Reduction0},
\eqref{eqn:Division3}, and
\eqref{eqn:Reduction3}.

Finally, we must look at the \textit{irregular} situation
with finitely many small values of the parameter \(r<L(b)\),
for fixed absolute coefficient \(b\ge 2\),
where the period shrinks to exceptional length \(\ell=2\).
While the initialization, \(0\)th reduction and \(1\)st division
remain the same as in Formulae
\eqref{eqn:Initialization},
\eqref{eqn:Reduction0},
\eqref{eqn:Division1},
the divergence sets in at the end of the \(1\)st reduction,
where the basis elements \(\vartheta^1_g\) and \(\vartheta^1_h\) are twisted. \\
\(1\)st \textbf{Reduction}: After execution of the steps from 1) to 6),
we get
\begin{equation}
\label{eqn:Reduction1Clone}
(1,\theta,\theta^2)\cdot\frac{1}{b}
\begin{pmatrix}
b & 0 & 0 \\
0 & 0 & b \\
0 & 1 & 0
\end{pmatrix}
=(1,\vartheta^1_g,\vartheta^1_h), \text{ the \(1\)st reduced Voronoi basis.}
\end{equation}
Thus, the second \textbf{relative} lattice minimum is
\(\vartheta^1_g=\theta^2/b\)
and its supplement is
\(\vartheta^1_h=b\theta/b=\theta\).
The second \textbf{absolute} lattice minimum is the product
\(\nu^2(1)=\vartheta^0_g\vartheta^1_g=\theta\cdot\theta^2/b
=\theta^3/b=1-3r\theta=\varepsilon_0\),
which is the desired \textbf{fundamental unit},
and the period length shrinks to \(\ell=2\). \\
\(2\)nd \textbf{Division}:
According to Formula
\eqref{eqn:Crucial},
\((\vartheta^1_g)^\prime(\vartheta^1_g)^{\prime\prime}
=(9r^2b^2+b\theta+3rb\theta^2)/b^2=(9r^2b+\theta+3r\theta^2)/b\)
and \(\mathrm{N}(\vartheta^1_g)=\mathrm{N}(\theta^2/b)=b^2/b^3=1/b\).
Thus, \(1/\vartheta^1_g=((9r^2b+\theta+3r\theta^2)/b)/(1/b)
=9r^2b+\theta+3r\theta^2\) and
\(\vartheta^1_h/\vartheta^1_g=\theta(9r^2b+\theta+3r\theta^2)
=3rb+\theta^2\),
and we get the \(2\)nd divided basis
\begin{equation}
\label{eqn:Division2Clone}
(1,\theta,\theta^2)\cdot
\begin{pmatrix}
1 & 3rb & 9r^2b \\
0 & 0 & 1 \\
0 & 1 & 3r
\end{pmatrix}
=(1,\vartheta^1_h/\vartheta^1_g,1/\vartheta^1_g).
\end{equation}
\(2\)nd \textbf{Reduction}: After execution of the steps from 1) to 6),
we get
\begin{equation}
\label{eqn:Reduction2Clone}
(1,\theta,\theta^2)\cdot
\begin{pmatrix}
1 & 0 & 1 \\
0 & 1 & -1 \\
0 & 0 & 1
\end{pmatrix}
=(1,\vartheta^2_g,\vartheta^2_h), \text{ the \(2\)nd reduced Voronoi basis.}
\end{equation}
This is the same lattice as in Formula
\eqref{eqn:Reduction0}
and the Voronoi algorithm can be terminated.
\end{proof}

\begin{remark}
\label{rmk:VoronoiPeriod}
For \(\ell\in\lbrace 1,3\rbrace\),
the reduction steps 1), 2) and 4) can always be skipped,
because the geometric conditions concerning the
hexagon of Zelling are satisfied a priori.
However, step 3) is required
in the \(1\)st and \(3\)rd reduction,
and step 5) shortens the longer basis vector (or both) by an integer
in the \(1\)st, \(2\)nd and \(3\)rd reduction.
For \(\ell=2\),
the twisting in step 1) is required in the \(2\)nd reduction.

We emphasize that there are also some typographical errors in
\cite[pp. 445--446]{WiZa1972}:
in step 4) the six pairs should begin with \((b,d)\) \((-b,-d)\), \(\ldots\),
in step 5) the first reduction should be \(\bar{\mu}=\mu-\lbrack\mu\rbrack\),
in step 6) the letter \(\rho\) is missing in the formula for \(\delta(\theta)\),
and the crucial selection of \(\theta_g\) and \(\theta_h\) is dependent on
\(\delta(\theta_g)<\delta(\theta_h)\) rather than \(\theta_g<\theta_h\).
In the final formula \(\epsilon_0=\prod_{i=1}^j\,\theta_g^{(i)}\),
the letter \(\theta\) is missing.

The relation which is responsible for the \textbf{shrinking of the period length}
from \(\ell=3\) to \(\ell=2\) is an exceptional behavior of the complements of
\(\vartheta_g^1=\theta^2/b\) and \(\vartheta_h^1=\theta\):
\((\vartheta_g^1)^\prime(\vartheta_g^1)^{\prime\prime}=(\theta^2/b)^\prime(\theta^2/b)^{\prime\prime}
=(9r^2b^2+b\theta+3rb\theta^2)/b^2\)
and \((\vartheta_h^1)^\prime(\vartheta_h^1)^{\prime\prime}=\theta^\prime\theta^{\prime\prime}=3rb+\theta^2\).
Since \(\theta^4=b\theta-3rb\theta^2\), the square of \(t:=3rb+\theta^2\)
is exactly \(t^2=9r^2b^2+b\theta+3rb\theta^2\), and the critical relation is
\(t^2/b^2<t\), which is equivalent to \(3rb+\theta^2<b^2\), respectively \(r<\frac{b}{3}-\frac{\theta^2}{3b}\).
For the interesting case \(b\ge 2\), we have \(0<\frac{\theta^2}{3b}<\frac{1}{6}<\frac{1}{3}\), and thus
we get the \textbf{sharp lower bound} \(L(b)\) as a ceiling integer:
\begin{equation}
\label{eqn:}
\begin{aligned}
r < L(b)=\lceil b/3\rceil & \Longleftrightarrow \ell=2 \text{ and the Voronoi chain is an \(\mathrm{M}1\)-chain}, \\
r\ge L(b)=\lceil b/3\rceil & \Longleftrightarrow \ell=3 \text{ and the Voronoi chain is an \(\mathrm{M}2\)-chain}.
\end{aligned}
\end{equation}
\end{remark}


\noindent
Lemma
\ref{lem:Monogenic}
demonstrates impressively that the concept of \textit{monogeneity}
has astonishing and powerful arithmetical applications
and is not a useless intellectual fancy theory.
We have stated the Lemma
in a form sufficiently general to admit a striking generalization
of item (1) in the Main Theorem
\ref{thm:Main}
to a number field \(L\) of odd prime degree \(p\ge 5\)
whose splitting field \(N\) has the maximal possible degree \(\lbrack N:\mathbb{Q}\rbrack=p!\)
since the Galois group \(\mathrm{Gal}(N/\mathbb{Q})\)
is the \textit{full symmetric group} \(S_p\).
Although the present article is mainly devoted to cubic number fields
we feel the desire to add a highlight concerning fields of higher degree
for which no general theory of \textit{differential principal factorizations}
is available yet,
let alone a general theory of \textit{multiplicity and multiplets}.


\begin{corollary}
\label{cor:Symmetric}
\textbf{(Absolute principal factors in \(S_p\)-fields.)}
Let \(p\ge 5\) be an arbitrary odd prime number bigger than \(3\).
A zero \(\theta\) of the \textbf{monogenic} polynomial
\(P(X)=X^p+\sigma\cdot prbX-b\in\mathbb{Z}\lbrack X\rbrack\)
with sign \(\sigma\in\lbrace -1,+1\rbrace\),
\textbf{squarefree} coefficient \(b\in\mathbb{N}\)
and parameter \(r\in\mathbb{N}\), such that
\(v_p(b^{p-1}+\sigma\cdot prb-1)=1\) and
\((\forall\,\ell\in\mathbb{P})\) \(\lbrack v_\ell(\sigma\cdot (p-1)^{p-1}r^pb+1)\ge 2\Longrightarrow\ell\mid prb\rbrack\),
generates a non-Galois field \(L=\mathbb{Q}(\theta)\) of degree \(p\) with maximal order
\(\mathcal{O}_L=\mathbb{Z}\lbrack\theta\rbrack\),
discriminant \(d_L=(-1)^{\frac{p(p-1)}{2}}p^pb^{p-1}\left(\sigma\cdot (p-1)^{p-1}r^pb+1\right)\),
and the following arithmetic properties.
For the examples in Table
\ref{tbl:Symmetric},
the Galois closure \(N\) of \(L\) possesses the \textbf{full symmetric group}
\(\mathrm{Gal}(N/\mathbb{Q})\simeq S_p\).
Generally, if the polynomial coefficient \(b=q\in\mathbb{P}\) is a prime number, then
\(\mathfrak{Q}=q\mathcal{O}_L+\theta\mathcal{O}_L\in\mathbb{P}_L\) is an ambiguous principal ideal,
that is, an \textbf{absolute differential principal factor} of \(L/\mathbb{Q}\).
\end{corollary}

\begin{proof}
According to
\cite[Lem. 1, p. 581]{LlNt1983},
the polynomial decomposition with \(p\) identical factors \(G(X)=X\) in
\(P(X)=X^p+\sigma\cdot prbX-b\equiv X^p\,(\mathrm{mod}\,q)\)
modulo each prime divisor \(q\in\mathbb{P}\) of \(b\)
establishes a prime factorization
\(q\mathcal{O}_L=\mathfrak{Q}^p\)
into the \(p\)-th power of the prime ideal
\(\mathfrak{Q}=q\mathcal{O}_L+\theta\mathcal{O}_L\in\mathbb{P}_L\) of \(\mathcal{O}_L\),
where trivially \(\theta=G(\theta)\),
that is, \(q\) is totally ramified in \(L\).
In the special case that \(b=q\) is itself prime, we can apply Lemma
\ref{lem:Monogenic},
which implies that
\(\mathfrak{Q}=q\mathcal{O}_L+\theta\mathcal{O}_L=\theta\mathcal{O}_L\)
is the ambiguous principal ideal generated by \(\theta\) in \(L\).
This is an \textbf{absolute principal factor} of \(L\)
\cite{Ma2019,Ma2021a,Ma2021b}.
\end{proof}


\renewcommand{\arraystretch}{1.1}
\begin{table}[ht]
\caption{Trinomials \(P(X)=X^p+aX-q\), \(a=\sigma\cdot prq\), with symmetric group \(S_p\)}
\label{tbl:Symmetric}
\begin{center}
\begin{tabular}{|r||c|r|r|r||r|l|c|r|}
\hline
 \(p\) & \(\sigma\) & \(r\) & \(q\) & \(a\)   & \(d_L\)                 & Factorization                               & DPF     & \(\#\mathfrak{G}\)  \\
\hline
 \(5\) &     \(+1\) & \(2\) & \(2\) &  \(20\) &       \(819\,250\,000\) & \(=2^4\cdot 5^5\cdot 5\cdot 29\cdot 113\)   & \((2)\) &  \(120\) \\
 \(5\) &     \(+1\) & \(1\) & \(3\) &  \(15\) &       \(194\,653\,125\) & \(=3^4\cdot 5^5\cdot 769\)                  & \((3)\) &  \(120\) \\
 \(5\) &     \(+1\) & \(1\) & \(7\) &  \(35\) &   \(13\,453\,103\,125\) & \(=7^4\cdot 5^5\cdot 11\cdot 163\)          & \((7)\) &  \(120\) \\
\hline
 \(5\) &     \(-1\) & \(1\) & \(2\) & \(-10\) &       \(-25\,550\,000\) & \(=-2^4\cdot 5^5\cdot 7\cdot 73\)           & \((2)\) &  \(120\) \\
 \(5\) &     \(-1\) & \(1\) & \(3\) & \(-15\) &      \(-194\,146\,875\) & \(=-3^4\cdot 5^5\cdot 13\cdot 59\)          & \((3)\) &  \(120\) \\
 \(5\) &     \(-1\) & \(2\) & \(7\) & \(-70\) & \(-430\,251\,696\,875\) & \(=-7^4\cdot 5^5\cdot 11\cdot 13\cdot 401\) & \((7)\) &  \(120\) \\
\hline
 \(7\) &     \(+1\) & \(1\) & \(2\) &  \(14\) &      \(-4\,918\,225\,149\,376\) & \(=-2^6\cdot 7^7\cdot 11\cdot 17\cdot 499\) & \((2)\) & \(5040\) \\
 \(7\) &     \(+1\) & \(1\) & \(3\) &  \(21\) &     \(-84\,032\,187\,331\,743\) & \(=-3^6\cdot 7^7\cdot 139\,969\)            & \((3)\) & \(5040\) \\
 \(7\) &     \(+1\) & \(1\) & \(5\) &  \(35\) & \(-3\,001\,827\,102\,859\,375\) & \(=-5^6\cdot 7^7\cdot 263\cdot 887\)        & \((5)\) & \(5040\) \\
\hline
 \(7\) &     \(-1\) & \(2\) & \(2\) & \(-28\) &     \(629\,526\,019\,949\,120\) & \(=2^6\cdot 7^7\cdot 5\cdot 691\cdot 3457\) & \((2)\) & \(5040\) \\
 \(7\) &     \(-1\) & \(1\) & \(3\) & \(-21\) &      \(84\,030\,986\,606\,049\) & \(=3^6\cdot 7^7\cdot 139\,967\)             & \((3)\) & \(5040\) \\
 \(7\) &     \(-1\) & \(1\) & \(5\) & \(-35\) &  \(3\,001\,801\,367\,140\,625\) & \(=5^6\cdot 7^7\cdot 233\,279\)             & \((5)\) & \(5040\) \\
\hline
\end{tabular}
\end{center}
\end{table}


\begin{example}
\label{exm:Symmetric}
In Table
\ref{tbl:Symmetric},
we give examples of quintic and septic number fields \(L=\mathbb{Q}(\theta)\),
generated by a zero \(\theta\) of a monogenic trinomial
\(P(X)=X^p+aX-b\in\mathbb{Z}\lbrack X\rbrack\)
with \(p\in\lbrace 5,7\rbrace\), \(a=\sigma\cdot prq\),
\(\sigma\in\lbrace -1,+1\rbrace\), \(r\in\mathbb{N}\), and
a \textit{prime} number \(b=q\in\mathbb{P}\).
The polynomials satisfy the assumptions of Corollary
\ref{cor:Symmetric},
and therefore the stem field \(L\) possesses an \textit{absolute} DPF,
\(\theta\mathcal{O}_L\in\mathbb{P}_L\), such that \((q)=q\mathcal{O}_L=(\theta\mathcal{O}_L)^p\),
and the splitting field \(N\) has a Galois group
\(\mathfrak{G}=\mathrm{Gal}(N/\mathbb{Q})\simeq S_p\) with order \(p!=\#\mathfrak{G}\).
Note that in contrast to the tables of section \S\
\ref{s:Experiments},
we are not able to use the concepts
conductor \(f\),
unique quadratic subfield \(K\) with discriminant \(d_K\) and \(p\)-class group \(\mathrm{Cl}_p(K)\),
logarithmic index \(w\) of subfield units,
and we do not compute the regulator in order to identify the relevant member of a multiplet,
since there are no definitions of the multiplicities \(m,m^\prime,m^{\prime\prime}\).
Only the discriminant \(d_L\), divisible by \(q^{p-1}p^p\), characterizes the ramification.
\end{example}


\section{Experimental verification of simply real cubic fields}
\label{s:Experiments}

\noindent
By means of the computational algebra system Magma
\cite{BCP1997,BCFS2022,MAGMA2022},
we have verified our theoretical statements experimentally.
In each of the following tables, we denote by
\(m\) the multiplicity of the \textit{homogeneous} multiplet \((L_1,\ldots,L_m)\) of cubic fields with conductor \(f\),
\(m^\prime\) the multiplicity of the \textit{heterogeneous} multiplet of \(3\)-ring class fields with conductors \(c\mid f\),
\(m^{\prime\prime}\) the multiplicity of the \textit{heterogeneous} multiplet of \(3\)-ray class fields with conductors \(c\mid f\).
The regulator \(\mathrm{Reg}=\log(\varepsilon_0^{-1})\) of \(L\)
identifies the member \(L_i\simeq L\) of the multiplet \((L_1,\ldots,L_m)\) uniquely.
Column \(d_K\) characterizes the quadratic fundamental discriminant by its residue class
\(\equiv d_K\,(\mathrm{mod}\,9)\).
We begin with \textit{simply real cubic} fields, \(\sigma=+1\), in Sections \S\S\
\ref{ss:Five}--\ref{ss:One}.
Since the type \(\beta\) of \(L\) and
the index \((U_N:U_0)=3^w\) of subfield units in the unit group of the Galois closure \(N\) of \(L\)
are constant with exponent \(w=1\), they are not listed.
The generating trinomial of \(L\) is always \(P(X)=X^3+aX-b\) with \(a=3rb\).
The primitive period length of the Voronoi algorithm is denoted by \(\ell\),
and the statement of Theorem
\ref{thm:VoronoiPeriod}
that an M2-chain occurs precisely for \(\ell=3\) is confirmed.


\subsection{Composite coefficients \(b=q_1\cdots q_k\)}
\label{ss:Five}

\noindent
The most extreme situations occur when \(b\) is selected as a \textit{primorial}
\(b=\prod\lbrace q\in\mathbb{P}\mid q\le U\rbrace\) with some upper bound \(U\),
for instance \(b=2\cdot 3\cdot 5\cdot 7\cdot 11=2310\) with \(U=11\).
It realizes the biggest multiplicities and principal factors with minimal theoretical effort.
The required experimental resources, however, are overwhelming.
For instance, the CPU-time for the computation of the six multiplets with \(m=48\) in Table
\ref{tbl:Five}
was more than two days and needed \(7\) Gigabytes of RAM,
since \(3280=(3^8-1)/2=m^{\prime\prime}\) ray class fields had to be constructed in the background,
for each multiplet.
According to Theorem
\ref{thm:VoronoiPeriod},
the sharp lower bound for \(\ell=3\) in Table
\ref{tbl:Five}
is \(r\ge L(b)=\lceil b/3\rceil=2310/3=770\),
far outside of the range of our investigations.

\begin{example}
\label{exm:Five}
In Table
\ref{tbl:Five},
both, the absolute coefficient
\(b=2310=3\cdot 2\cdot 5\cdot 7\cdot 11\)
and the conductor
\(f=6930=3^2\cdot 2\cdot 5\cdot 7\cdot 11\)
have five prime divisors.
With respect to
\cite[Thm. 3.2, p. 2215]{Ma2014},
\(f=6930=3^2\cdot 770\) is a \textit{free irregular} conductor with \(\delta_3(f)=0\), \(\omega=1\), \(\tau=5\)
and multiplicity
\(m=m_3(d_K,f)=3^{\varrho+\omega}\cdot 2^{\tau-1}=3^{0+1}\cdot 2^{5-1}=3\cdot 16=48\).
The norm \(\mathrm{N}_{L/\mathbb{Q}}(\theta)=b=2310\) of the absolute DPF \(\theta\) of \(L/\mathbb{Q}\) and
the fundamental unit \(\varepsilon_0=1-3r\theta<1\) of \(L\) verify Theorem
\ref{thm:TrinomialPrinciple}
with \(v=3r\), \(d=3\) and \(t=1\).
Here, the \(3\)-class group \(\mathrm{Cl}_3(K)=1\) of \(K\) is always trivial
with rank \(\varrho=\varrho_3(K)=0\).
We point out that \(m^\prime=(3^{\varrho+\tau+\omega-\delta}-1)/2=(3^{0+5+1-0}-1)/2=364\)
\cite{Ma2014}.
\end{example}

\renewcommand{\arraystretch}{1.1}
\begin{table}[ht]
\caption{Irregular conductor \(f\) with \(\tau=5\) and most extensive multiplets}
\label{tbl:Five}
\begin{center}
\begin{tabular}{|r|r|r||c|c|c|r|c|c|c||r|r|r|c|}
\hline
 \(b\)    & \(r\) & \(a\)     & \(f\)    & \(d_K\) & \(\varrho\) & \(\varepsilon_0\) & \(\mathrm{Reg}\) & \(\ell\) & Chain & \(m^{\prime\prime}\) & \(m^\prime\) & \(m\)  & Multiplet     \\
\hline
 \(2310\) & \(1\) &  \(6930\) & \(6930\) &  \(-3\) &       \(0\) &     \(1-3\theta\) &       \(11.041\) &    \(2\) & M1 &             \(3280\) &      \(364\) & \(48\) & \((\beta^{48})\) \\
 \(2310\) & \(2\) & \(13860\) & \(6930\) &  \(-3\) &       \(0\) &     \(1-6\theta\) &       \(13.120\) &    \(2\) & M1 &             \(3280\) &      \(364\) & \(48\) & \((\beta^{48})\) \\
 \(2310\) & \(3\) & \(20790\) & \(6930\) &  \(-3\) &       \(0\) &     \(1-9\theta\) &       \(14.337\) &    \(2\) & M1 &             \(3280\) &      \(364\) & \(48\) & \((\beta^{48})\) \\
 \(2310\) & \(5\) & \(34650\) & \(6930\) &  \(-3\) &       \(0\) &    \(1-15\theta\) &       \(15.869\) &    \(2\) & M1 &             \(3280\) &      \(364\) & \(48\) & \((\beta^{48})\) \\
 \(2310\) & \(6\) & \(41580\) & \(6930\) &  \(-3\) &       \(0\) &    \(1-18\theta\) &       \(16.416\) &    \(2\) & M1 &             \(3280\) &      \(364\) & \(48\) & \((\beta^{48})\) \\
 \(2310\) & \(7\) & \(48510\) & \(6930\) &  \(-3\) &       \(0\) &    \(1-21\theta\) &       \(16.879\) &    \(2\) & M1 &             \(3280\) &      \(364\) & \(48\) & \((\beta^{48})\) \\
\hline 
\end{tabular}
\end{center}
\end{table}


\noindent
In Example
\ref{exm:Five},
the divisibility of \(b\) by \(3\) is annoying,
since it enforces an irregular conductor \(f\), by Theorem
\ref{thm:Preparation}.
Therefore, we choose a modified primorial \(b=2\cdot 5\cdot 7\cdot 11=770\),
liberated from the contribution by the critical prime divisor \(3\),
and thus with four prime divisors only, in Example
\ref{exm:FourFive}.

\begin{example}
\label{exm:FourFive}
In Table
\ref{tbl:FourFive},
there appear two conductors \(6930\) and \(2310\),
and the columns are arranged in the same way as in Table
\ref{tbl:Five}.
The conductor
\(f=2310=3\cdot 2\cdot 5\cdot 7\cdot 11\)
has the same five prime divisors as \(6930\) but is squarefree.
The norm \(\mathrm{N}_{L/\mathbb{Q}}(\theta)=b=770\) of the absolute DPF \(\theta\) of \(L/\mathbb{Q}\) and
the fundamental unit \(\varepsilon_0=1-3r\theta\) of \(L\) verify Theorem
\ref{thm:TrinomialPrinciple}.
When the quadratic fundamental discriminant \(d_K\) is coprime to \(3\),
the conductor must be \(f=9b=6930\).
When \(d_K\) is a multiple of \(3\), the conductor must be \(f=3b=2310\).
In particular, \(f\) is \textit{never irregular}, according to Theorem
\ref{thm:Preparation}.
For the first five cases \(r\in\lbrace 1,4,7,10\rbrace\) of \(f=6930\), 
and for the third and fifth case \(r\in\lbrace 9,15\rbrace\) of \(f=2310\),
we have \(\varrho=0\) and thus
\(m=3^{\varrho+\omega}\cdot 2^{\tau-1}=3^{0+0}\cdot 2^{5-1}=16\).
However, positive \(3\)-class rank \(\varrho=1\) of \(K\) appears occasionally,
and then we have to use the \textit{restrictive multiplicity formula}
\cite[Thm. 3.3, p. 2217]{Ma2014},
\(m=3^{\varrho+\omega}\cdot 2^{u}\cdot\frac{1}{3}\lbrack 2^{v-1}-(-1)^{v-1}\rbrack\),
where the beginning of the sequence of restrictive contributions
\((\frac{1}{3}\lbrack 2^{v-1}-(-1)^{v-1}\rbrack)_{v\ge 0}\)
is given by \((\frac{1}{2},0,1,1,3,5,\ldots)\),
and \(\tau=u+v\) expresses the partition into
\textit{free} and \textit{restrictive} prime divisors of \(f\).
Since \(\tau=5\), the free formula arises as degenerate case \((u,v)=(5,0)\),
with multiplicity \(m=3^{1+0}\cdot 2^{5-1}=3\cdot 16=48\), which does not occur.
So the experimental result for the multiplicity \(m\) decides that
a \textit{restriction} \(\delta_3(f)=1\) of \textit{first order} occurs.
Since \((u,v)=(4,1)\) is impossible, we must calculate \(m\) for the other pairs:
\begin{itemize}
\item
\((u,v)=(3,2)\) yields \(m=3^{1+0}\cdot 2^{3}\cdot 1=24\),
\item
\((u,v)=(2,3)\) yields \(m=3^{1+0}\cdot 2^{2}\cdot 1=12\),
\item
\((u,v)=(1,4)\) yields \(m=3^{1+0}\cdot 2\cdot 3=18\),
\item
\((u,v)=(0,5)\) yields \(m=3^{1+0}\cdot 1\cdot 5=15\).
\end{itemize} 
While \(m=15\) does not arise, we must have
\((u,v)=(1,4)\) for \(r=13\),
\((u,v)=(2,3)\) for \(r\in\lbrace 3,6\rbrace\), and
\((u,v)=(3,2)\) for \(r=12\).
Note that \(m^\prime=(3^{\varrho+5+0-\delta}-1)/2=121\) is constant, since \(\varrho=\delta\).
\end{example}

\renewcommand{\arraystretch}{1.1}
\begin{table}[ht]
\caption{Regular conductors \(f\) with \(\tau=5\) and smaller multiplets}
\label{tbl:FourFive}
\begin{center}
\begin{tabular}{|r|r|r||c|c|c|r|r|c|c||r|r|r|c|c|}
\hline
 \(b\)   & \(r\)  & \(a\)     & \(f\)    & \(d_K\) & \(\varrho\) & \(\varepsilon_0\) & \(\mathrm{Reg}\) & \(\ell\) & Chain & \(m^{\prime\prime}\) & \(m^\prime\) & \(m\)  & Multiplet     & \(v\) \\
\hline
 \(770\) &  \(1\) &  \(2310\) & \(6930\) &  \(-1\) &       \(0\) &     \(1-3\theta\) &        \(9.942\) &    \(2\) & M1 &             \(1093\) &      \(121\) & \(16\) & \((\beta^{16})\) & \(0\) \\
 \(770\) &  \(4\) &  \(9240\) & \(6930\) &  \(-7\) &       \(0\) &    \(1-12\theta\) &       \(14.101\) &    \(2\) & M1 &             \(1093\) &      \(121\) & \(16\) & \((\beta^{16})\) & \(0\) \\
 \(770\) &  \(7\) & \(16170\) & \(6930\) &  \(-4\) &       \(0\) &    \(1-21\theta\) &       \(15.780\) &    \(2\) & M1 &             \(1093\) &      \(121\) & \(16\) & \((\beta^{16})\) & \(0\) \\
 \(770\) & \(10\) & \(23100\) & \(6930\) &  \(-1\) &       \(0\) &    \(1-30\theta\) &       \(16.850\) &    \(2\) & M1 &             \(1093\) &      \(121\) & \(16\) & \((\beta^{16})\) & \(0\) \\
 \(770\) & \(13\) & \(30030\) & \(6930\) &  \(-7\) &       \(1\) &    \(1-39\theta\) &       \(17.637\) &    \(2\) & M1 &             \(1093\) &      \(121\) & \(18\) & \((\beta^{18})\) & \(4\) \\
\hline 
 \(770\) &  \(3\) &  \(6930\) & \(2310\) &  \(-3\) &       \(1\) &     \(1-9\theta\) &       \(13.238\) &    \(2\) & M1 &              \(364\) &      \(121\) & \(12\) & \((\beta^{12})\) & \(3\) \\
 \(770\) &  \(6\) & \(13860\) & \(2310\) &  \(-3\) &       \(1\) &    \(1-18\theta\) &       \(15.318\) &    \(2\) & M1 &              \(364\) &      \(121\) & \(12\) & \((\beta^{12})\) & \(3\) \\
 \(770\) &  \(9\) & \(20790\) & \(2310\) &  \(-3\) &       \(0\) &    \(1-27\theta\) &       \(16.534\) &    \(2\) & M1 &              \(364\) &      \(121\) & \(16\) & \((\beta^{16})\) & \(0\) \\
 \(770\) & \(12\) & \(27720\) & \(2310\) &  \(-3\) &       \(1\) &    \(1-36\theta\) &       \(17.397\) &    \(2\) & M1 &              \(364\) &      \(121\) & \(24\) & \((\beta^{24})\) & \(2\) \\
 \(770\) & \(15\) & \(34650\) & \(2310\) &  \(-3\) &       \(0\) &    \(1-45\theta\) &       \(18.066\) &    \(2\) & M1 &              \(364\) &      \(121\) & \(16\) & \((\beta^{16})\) & \(0\) \\
\hline
\end{tabular}
\end{center}
\end{table}


\renewcommand{\arraystretch}{1.1}
\begin{table}[ht]
\caption{Regular conductor \(f\) with \(\tau=4\) and modest multiplets}
\label{tbl:Four}
\begin{center}
\begin{tabular}{|r|r|r||c|c|c|r|r|c|c||r|r|r|c|c|}
\hline
 \(b\)  & \(r\)  & \(a\)    & \(f\)   & \(d_K\) & \(\varrho\) & \(\varepsilon_0\)     & \(\mathrm{Reg}\) & \(\ell\) & Chain & \(m^{\prime\prime}\) & \(m^\prime\) & \(m\)  & Multiplet     & \(v\) \\
\hline
 \(70\) &  \(1\) &  \(210\) & \(210\) &  \(-6\) &       \(1\) &         \(1-3\theta\) &        \(7.546\) &    \(2\) & M1 &              \(121\) &       \(40\) &  \(9\) &  \((\beta^{9})\) & \(4\) \\
 \(70\) &  \(3\) &  \(630\) & \(210\) &  \(-3\) &       \(0\) &         \(1-9\theta\) &       \(10.840\) &    \(2\) & M1 &              \(121\) &       \(40\) &  \(8\) &  \((\beta^{8})\) & \(0\) \\
 \(70\) &  \(4\) &  \(840\) & \(210\) &  \(-6\) &       \(0\) &        \(1-12\theta\) &       \(11.703\) &    \(2\) & M1 &              \(121\) &       \(40\) &  \(8\) &  \((\beta^{8})\) & \(0\) \\
 \(70\) &  \(6\) & \(1260\) & \(210\) &  \(-3\) &       \(1\) &        \(1-18\theta\) &       \(12.920\) &    \(2\) & M1 &              \(121\) &       \(40\) &  \(9\) & \((\alpha_2^{4},\beta^{5})\) & \(4\) \\
 \(70\) &  \(7\) & \(1470\) & \(210\) &  \(-6\) &       \(0\) &        \(1-21\theta\) &       \(13.382\) &    \(2\) & M1 &              \(121\) &       \(40\) &  \(8\) &  \((\beta^{8})\) & \(0\) \\
 \(70\) &  \(9\) & \(1890\) & \(210\) &  \(-3\) &       \(0\) &        \(1-27\theta\) &       \(14.136\) &    \(2\) & M1 &              \(121\) &       \(40\) &  \(8\) &  \((\alpha_2^{4},\beta^{4})\) & \(0\) \\
 \(70\) & \(10\) & \(2100\) & \(210\) &  \(-6\) &       \(1\) &        \(1-30\theta\) &       \(14.452\) &    \(2\) & M1 &              \(121\) &       \(40\) &  \(9\) &  \((\beta^{9})\) & \(4\) \\
 \(70\) & \(12\) & \(2520\) & \(210\) &  \(-3\) &       \(1\) &        \(1-36\theta\) &       \(14.999\) &    \(2\) & M1 &              \(121\) &       \(40\) & \(12\) & \((\beta^{12})\) & \(2\) \\
 \(70\) & \(13\) & \(2730\) & \(210\) &  \(-6\) &       \(0\) &        \(1-39\theta\) &       \(15.239\) &    \(2\) & M1 &              \(121\) &       \(40\) &  \(8\) &  \((\beta^{8})\) & \(0\) \\
 \(70\) & \(15\) & \(3150\) & \(210\) &  \(-3\) &       \(0\) &        \(1-45\theta\) &       \(15.668\) &    \(2\) & M1 &              \(121\) &       \(40\) &  \(8\) &  \((\beta^{8})\) & \(0\) \\
 \(70\) & \(16\) & \(3360\) & \(210\) &  \(-6\) &       \(0\) & \((1-48\theta)^{-1}\) &       \(15.862\) &    \(2\) & M1 &              \(121\) &       \(40\) &  \(8\) &  \((\beta^{8})\) & \(0\) \\
 \(70\) & \(18\) & \(3780\) & \(210\) &  \(-3\) &       \(1\) &        \(1-54\theta\) &       \(16.215\) &    \(2\) & M1 &              \(121\) &       \(40\) & \(12\) & \((\beta^{12})\) & \(2\) \\
 \(70\) & \(19\) & \(3990\) & \(210\) &  \(-6\) &       \(1\) & \((1-57\theta)^{-1}\) &       \(16.378\) &    \(2\) & M1 &              \(121\) &       \(40\) &  \(9\) &  \((\beta^{9})\) & \(4\) \\
 \(70\) & \(21\) & \(4410\) & \(210\) &  \(-3\) &       \(0\) & \((1-63\theta)^{-1}\) &       \(16.678\) &    \(2\) & M1 &              \(121\) &       \(40\) &  \(8\) &  \((\beta^{8})\) & \(0\) \\
 \(70\) & \(22\) & \(4620\) & \(210\) &  \(-6\) &       \(0\) &        \(1-66\theta\) &       \(16.817\) &    \(2\) & M1 &              \(121\) &       \(40\) &  \(8\) &  \((\beta^{8})\) & \(0\) \\
 \(70\) & \(\mathbf{24}\) & \(5040\) & \(210\) & \(-3\) & \(0\) &      \(1-72\theta\) &       \(17.078\) & \(\mathbf{3}\) & M2 &        \(121\) &       \(40\) &  \(8\) &  \((\beta^{8})\) & \(0\) \\
 \(70\) & \(25\) & \(5250\) & \(210\) &  \(-6\) &       \(0\) &        \(1-75\theta\) &       \(17.201\) &    \(3\) & M2 &              \(121\) &       \(40\) &  \(8\) &  \((\beta^{8})\) & \(0\) \\
 \(70\) & \(27\) & \(5670\) & \(210\) &  \(-3\) &       \(0\) &        \(1-81\theta\) &       \(17.432\) &    \(3\) & M2 &              \(121\) &       \(40\) &  \(8\) &  \((\beta^{8})\) & \(0\) \\
 \(70\) & \(28\) & \(5880\) & \(210\) &  \(-6\) &       \(1\) &        \(1-84\theta\) &       \(17.541\) &    \(3\) & M2 &              \(121\) &       \(40\) & \(12\) & \((\beta^{12})\) & \(2\) \\
 \(70\) & \(30\) & \(6300\) & \(210\) &  \(-3\) &       \(1\) & \((1-90\theta)^{-1}\) &       \(17.748\) &    \(3\) & M2 &              \(121\) &       \(40\) &  \(6\) &  \((\beta^{6})\) & \(3\) \\
\hline
\end{tabular}
\end{center}
\end{table}

\noindent
The computations for Tables
\ref{tbl:Five}--\ref{tbl:FourFive}
squeezed the computational algebra system Magma
\cite{MAGMA2022}
close to its limits.
They had to be terminated due to unmanageable
CPU times and RAM storage requirements.
Thus, we take a smaller modified primorial for the absolute coefficient
\(b=2\cdot 5\cdot 7=70\),
again abstaining from the critical prime divisor \(3\).

\begin{example}
\label{exm:Four}
In Table
\ref{tbl:Four},
the conductor \(f=210\) is constant,
and the columns are arranged in the same way as in Table
\ref{tbl:Five}.
The norm \(\mathrm{N}_{L/\mathbb{Q}}(\theta)=b=70\) of the absolute DPF \(\theta\) of \(L/\mathbb{Q}\)
and the fundamental unit \(\varepsilon_0=1-3r\theta\) of \(L\) verify Theorem
\ref{thm:TrinomialPrinciple}.
The sharp bound where M2-chains set in is \(r=\lceil 70/3\rceil=24\).
For \(r\in\lbrace 16,19,21,30\rbrace\),
Magma computes the \textit{inverse} fundamental unit \(\varepsilon_0^{-1}>1\).
Since \(\tau=4\), the free formula arises as degenerated case \((u,v)=(4,0)\),
with multiplicity \(m=3^{0+0}\cdot 2^{4-1}=8\), if \(\varrho=0\).
For \(\varrho=1\), however,
the multiplicity is \(m=3^{1+0}\cdot 2^{4-1}=3\cdot 8=24\), which does not occur.
As before, the experimental result for the multiplicity \(m\) decides that
a \textit{restriction} \(\delta_3(f)=1\) of \textit{first order} occurs.
Since \((u,v)=(3,1)\) implies \(m=0\), we must calculate \(m\) for the other pairs:
\begin{itemize}
\item
\((u,v)=(2,2)\) yields \(m=3^{1+0}\cdot 2^{2}\cdot 1=12\),
\item
\((u,v)=(1,3)\) yields \(m=3^{1+0}\cdot 2\cdot 1=6\),
\item
\((u,v)=(0,4)\) yields \(m=3^{1+0}\cdot 1\cdot 3=9\).
\end{itemize} 
Actually, each case arises, and we must have
\((u,v)=(0,4)\) for \(r\in\lbrace 1,6,10,19\rbrace\),
\((u,v)=(1,3)\) for \(r=30\), and
\((u,v)=(2,2)\) for \(r\in\lbrace 12,18,28\rbrace\).
A \textit{relative} DPF occurs in four members with type \(\alpha_2\)
of the multiplets for \(r\in\lbrace 6,9\rbrace\),
enabled by the prime divisor \(7\) of \(b\),
which splits in \(K\).
As mentioned,  Magma occasionally computes the inverse of \(1-3r\theta\)
as fundamental unit, according to Formula
\eqref{eqn:Inverse},
\((1-3r\theta)^{-1}=(1+27r^3b)+3r\theta+9r^2\theta^2\).
For instance, we have 
\((1-90\theta)^{-1}=51\,030\,001+90\theta+8100\theta^2\) for \(r=30\).
\end{example}


\renewcommand{\arraystretch}{1.1}
\begin{table}[ht]
\caption{Regular conductors \(f\) with \(\tau=3\) and small multiplets}
\label{tbl:Three}
\begin{center}
\begin{tabular}{|r|r|r||c|c|c|r|r|c|c||r|r|r|c|c|}
\hline
 \(b\)  & \(r\)  & \(a\)   & \(f\)   & \(d_K\) & \(\varrho\) & \(\varepsilon_0\)     & \(\mathrm{Reg}\) & \(\ell\) & Chain & \(m^{\prime\prime}\) & \(m^\prime\) & \(m\)  & Multiplet     & \(v\) \\
\hline
 \(14\) &  \(1\) &  \(42\) & \(126\) &  \(-1\) &       \(0\) &         \(1-3\theta\) &        \(5.943\) &    \(2\) & M1 &              \(121\) &       \(13\) &  \(4\) &  \((\beta^{4})\) & \(0\) \\
 \(14\) &  \(4\) & \(168\) & \(126\) &  \(-7\) &       \(0\) &        \(1-12\theta\) &       \(10.094\) &    \(2\) & M1 &              \(121\) &       \(13\) &  \(4\) &  \((\beta^{4})\) & \(0\) \\
\hline
 \(14\) &  \(3\) & \(126\) &  \(42\) &  \(-3\) &       \(0\) &         \(1-9\theta\) &        \(9.231\) &    \(2\) & M1 &               \(40\) &       \(13\) &  \(4\) &  \((\beta^{4})\) & \(0\) \\
 \(14\) & \(\mathbf{6}\) & \(252\) & \(42\) & \(-3\) & \(0\) &        \(1-18\theta\) &       \(11.310\) & \(\mathbf{3}\) & M2 &         \(40\) &       \(13\) &  \(4\) &  \((\beta^{4})\) & \(0\) \\
\hline
\end{tabular}
\end{center}
\end{table}

\begin{example}
\label{exm:Three}
In Table
\ref{tbl:Three},
two conductors \(f=42=3\cdot 2\cdot 7\) and \(f=126=3^2\cdot 2\cdot 7\) appear.
The sharp bound where \(\mathrm{M}2\)-chains set in
is \(r\ge L(b)=\lceil b/3\rceil=\lceil 14/3\rceil=15/3=5\),
in fact \(r\ge 6\).
Since \(\varrho=0\), the multiplicity is \(m=3^{0+0}\cdot 2^{3-1}=2^2=4\).
\end{example}


\subsection{Prime coefficient \(b=q\)}
\label{ss:Two}

\noindent
In this case, the ambiguous principal ideal
\(\theta\mathcal{O}_L\) with norm \(\mathrm{N}_{L/\mathbb{Q}}(\theta)=q\in\mathbb{P}\)
in our Main Theorem
\ref{thm:Main}
is a prime ideal \(\mathfrak{Q}\in\mathbb{P}_L\)
such that \(q\mathcal{O}_L=\mathfrak{Q}^3\).

\renewcommand{\arraystretch}{1.1}
\begin{table}[ht]
\caption{Regular conductor \(f\) with \(\tau=2\) and small multiplets}
\label{tbl:Two}
\begin{center}
\begin{tabular}{|c|r|r||c|c|c|r|r|c|c|c||r|r|r|c|c|}
\hline
 \(b\) & \(r\)  & \(a\)   & \(f\)  & \(d_K\) & \(\varrho\) & \(\varepsilon_0\)     & \(\mathrm{Reg}\) & \(\mathrm{Cl}(L)\) & \(\ell\) & Chain & \(m^{\prime\prime}\) & \(m^\prime\) & \(m\) & Multiplet      & \(v\) \\
\hline
 \(7\) &  \(1\) &  \(21\) & \(21\) &  \(-6\) &       \(1\) &         \(1-3\theta\) &        \(5.257\) & \((3)\)            &    \(2\) & M1 &               \(13\) &        \(4\) &  \(3\) &  \((\beta^{3})\) & \(2\) \\
 \(7\) & \(\mathbf{3}\) & \(63\) & \(21\) & \(-3\) & \(0\) &  \((1-9\theta)^{-1}\) &        \(8.538\) & \((3,3)\)          & \(\mathbf{3}\) & M2 &         \(13\) &        \(4\) &  \(2\) &  \((\beta^{2})\) & \(0\) \\
 \(7\) &  \(4\) &  \(84\) & \(21\) &  \(-6\) &       \(1\) & \((1-12\theta)^{-1}\) &        \(9.401\) & \((3,12)\)         &    \(3\) & M2 &               \(13\) &        \(4\) &  \(3\) &  \((\alpha_2,\beta^{2})\) & \(2\) \\
 \(7\) &  \(6\) & \(126\) & \(21\) &  \(-3\) &       \(0\) & \((1-18\theta)^{-1}\) &       \(10.617\) & \((3,24)\)         &    \(3\) & M2 &               \(13\) &        \(4\) &  \(2\) &  \((\beta^{2})\) & \(0\) \\
 \(7\) &  \(7\) & \(147\) & \(21\) &  \(-6\) &       \(0\) & \((1-21\theta)^{-1}\) &       \(11.080\) & \((3,12)\)         &    \(3\) & M2 &               \(13\) &        \(4\) &  \(2\) &  \((\beta^{2})\) & \(0\) \\
 \(7\) &  \(9\) & \(189\) & \(21\) &  \(-3\) &       \(1\) & \((1-27\theta)^{-1}\) &       \(11.833\) & \((2,12)\)         &    \(3\) & M2 &               \(13\) &        \(4\) &  \(3\) &  \((\beta^{3})\) & \(2\) \\
 \(7\) & \(10\) & \(210\) & \(21\) &  \(-6\) &       \(1\) & \((1-30\theta)^{-1}\) &       \(12.150\) & \((132)\)          &    \(3\) & M2 &               \(13\) &        \(4\) &  \(3\) &  \((\beta^{3})\) & \(2\) \\
 \(7\) & \(12\) & \(252\) & \(21\) &  \(-3\) &       \(0\) & \((1-36\theta)^{-1}\) &       \(12.696\) & \((6,24)\)         &    \(3\) & M2 &               \(13\) &        \(4\) &  \(2\) &  \((\beta^{2})\) & \(0\) \\
 \(7\) & \(13\) & \(273\) & \(21\) &  \(-6\) &       \(0\) & \((1-39\theta)^{-1}\) &       \(12.937\) & \((48)\)           &    \(3\) & M2 &               \(13\) &        \(4\) &  \(2\) &  \((\beta^{2})\) & \(0\) \\
 \(7\) & \(16\) & \(336\) & \(21\) &  \(-6\) &       \(1\) & \((1-48\theta)^{-1}\) &       \(13.560\) & \((168)\)          &    \(3\) & M2 &               \(13\) &        \(4\) &  \(3\) &  \((\beta^{3})\) & \(2\) \\
 \(7\) & \(18\) & \(378\) & \(21\) &  \(-3\) &       \(0\) & \((1-54\theta)^{-1}\) &       \(13.913\) & \((156)\)          &    \(3\) & M2 &               \(13\) &        \(4\) &  \(2\) &  \((\beta^{2})\) & \(0\) \\
 \(7\) & \(19\) & \(399\) & \(21\) &  \(-6\) &       \(1\) & \((1-57\theta)^{-1}\) &       \(14.075\) & \((3,30)\)         &    \(3\) & M2 &               \(40\) &       \(13\) &  \(6\) &  \((\alpha_1,\beta^{5})\) & \(0\) \\
 \(7\) & \(21\) & \(441\) & \(21\) &  \(-3\) &       \(0\) & \((1-63\theta)^{-1}\) &       \(14.375\) & \((132)\)          &    \(3\) & M2 &               \(13\) &        \(4\) &  \(2\) &  \((\beta^{2})\) & \(0\) \\
 \(7\) & \(22\) & \(462\) & \(21\) &  \(-6\) &       \(1\) & \((1-66\theta)^{-1}\) &       \(14.515\) & \((3,117)\)        &    \(3\) & M2 &               \(40\) &       \(13\) &  \(6\) &  \((\alpha_1^5,\beta)\) & \(0\) \\
\hline
\end{tabular}
\end{center}
\end{table}

\begin{example}
\label{exm:Two}
In Table
\ref{tbl:Two},
the conductor \(f=21=3\cdot 7\) is constant.
The sharp bound where \(\mathrm{M}2\)-chains set in
is \(r\ge L(b)=\lceil b/3\rceil=\lceil 7/3\rceil=9/3=3\).
The multiplicity is \(m=2^{\tau-1}=2\) for \(\varrho=0\),
and either \(m=3^1\cdot 2^0\cdot\frac{1}{3}\lbrack 2^{2-1}-(-1)^{2-1}\rbrack=3\) with \((u,v)=(0,2)\)
or \(m=3^1\cdot 2^2\cdot\frac{1}{3}\lbrack 2^{0-1}-(-1)^{0-1}\rbrack=6\) with \((u,v)=(2,0)\) for \(\varrho=1\).
For \(r=22\), the field \(L=\mathbb{Q}(\theta)\) is \textbf{unique} with absolute DPF.
\end{example}


\renewcommand{\arraystretch}{1.1}
\begin{table}[ht]
\caption{Regular conductors \(f\) with \(\tau=1\) and tiny multiplets}
\label{tbl:One}
\begin{center}
\begin{tabular}{|r|r|r||c|c|c|r|r|c|c|c||r|r|r|c|c|}
\hline
 \(b\) & \(r\)  & \(a\)   & \(f\) & \(d_K\) & \(\varrho\) & \(\varepsilon_0\) & \(\mathrm{Reg}\) & \(\mathrm{Cl}(L)\) & \(\ell\) & Chain & \(m^{\prime\prime}\) & \(m^\prime\) & \(m\) & Multiplet        & \(v\) \\
\hline
 \(1\) &  \(1\) &   \(3\) & \(3\) &  \(-6\) &       \(0\) &        \(\theta\) &        \(1.133\) & \(1\)              &    \(1\) & M0 &                \(1\) &        \(1\) & \(1\) &      \((\beta)\) & \(0\) \\
\hline
 \(1\) &  \(2\) &   \(6\) & \(9\) &  \(-2\) &       \(0\) &        \(\theta\) &        \(1.796\) & \((3)\)            &    \(1\) & M0 &                \(4\) &        \(1\) & \(1\) &   \((\alpha_2)\) & \(0\) \\
 \(1\) &  \(4\) &  \(12\) & \(3\) &  \(-6\) &       \(1\) &        \(\theta\) &        \(2.485\) & \((6)\)            &    \(1\) & M0 &                \(4\) &        \(4\) & \(3\) & \((\alpha_1^3)\) & \(0\) \\
 \(1\) &  \(5\) &  \(15\) & \(9\) &  \(-5\) &       \(0\) &        \(\theta\) &        \(2.708\) & \((3)\)            &    \(1\) & M0 &                \(4\) &        \(1\) & \(1\) &   \((\alpha_2)\) & \(0\) \\
 \(1\) &  \(7\) &  \(21\) & \(3\) &  \(-6\) &       \(1\) &        \(\theta\) &        \(3.045\) & \((6)\)            &    \(1\) & M0 &                \(4\) &        \(4\) & \(3\) & \((\alpha_1^3)\) & \(0\) \\
 \(1\) &  \(8\) &  \(24\) & \(9\) &  \(-8\) &       \(0\) &        \(\theta\) &        \(3.178\) & \((18)\)           &    \(1\) & M0 &                \(4\) &        \(1\) & \(1\) &   \((\alpha_2)\) & \(0\) \\
 \(1\) & \(10\) &  \(30\) & \(3\) &  \(-6\) &       \(1\) &        \(\theta\) &        \(3.401\) & \((15)\)           &    \(1\) & M0 &                \(4\) &        \(4\) & \(3\) & \((\alpha_1^3)\) & \(0\) \\
 \(1\) & \(13\) &  \(39\) & \(3\) &  \(-6\) &       \(1\) &        \(\theta\) &        \(3.664\) & \((15)\)           &    \(1\) & M0 &                \(4\) &        \(4\) & \(3\) & \((\alpha_1^3)\) & \(0\) \\
 \(1\) & \(14\) &  \(42\) & \(9\) &  \(-5\) &       \(0\) &        \(\theta\) &        \(3.738\) & \((39)\)           &    \(1\) & M0 &                \(4\) &        \(1\) & \(1\) &   \((\alpha_2)\) & \(0\) \\
 \(1\) & \(16\) &  \(48\) & \(3\) &  \(-6\) &       \(1\) &        \(\theta\) &        \(3.871\) & \((48)\)           &    \(1\) & M0 &                \(4\) &        \(4\) & \(3\) & \((\alpha_1^3)\) & \(0\) \\
 \(1\) & \(17\) &  \(51\) & \(9\) &  \(-8\) &       \(1\) &        \(\theta\) &        \(3.932\) & \((3,3)\)          &    \(1\) & M0 &               \(13\) &        \(4\) & \(3\) & \((\alpha_1,\alpha_2^2)\) & \(0\) \\
 \(1\) & \(19\) &  \(57\) & \(3\) &  \(-6\) &       \(1\) &        \(\theta\) &        \(4.043\) & \((12)\)           &    \(1\) & M0 &                \(4\) &        \(4\) & \(3\) & \((\alpha_1^3)\) & \(0\) \\
 \(1\) & \(20\) &  \(60\) & \(9\) &  \(-2\) &       \(1\) &        \(\theta\) &        \(3.932\) & \((3,12)\)         &    \(1\) & M0 &               \(13\) &        \(4\) & \(3\) & \((\alpha_2^3)\) & \(0\) \\
 \(1\) & \(22\) &  \(66\) & \(3\) &  \(-6\) &       \(1\) &        \(\theta\) &        \(4.190\) & \((42)\)           &    \(1\) & M0 &                \(4\) &        \(4\) & \(3\) & \((\alpha_1^3)\) & \(0\) \\
 \(1\) & \(23\) &  \(69\) & \(9\) &  \(-5\) &       \(0\) &        \(\theta\) &        \(4.234\) & \((24)\)           &    \(1\) & M0 &                \(4\) &        \(1\) & \(1\) &   \((\alpha_2)\) & \(0\) \\
 \(1\) & \(25\) &  \(75\) & \(3\) &  \(-6\) &       \(2\) &        \(\theta\) &        \(4.317\) & \((6,6)\)          &    \(1\) & M0 &               \(13\) &       \(13\) & \(9\) & \((\alpha_1^9)\) & \(0\) \\
 \(1\) & \(26\) &  \(78\) & \(9\) &  \(-8\) &       \(0\) &        \(\theta\) &        \(4.357\) & \((2,36)\)         &    \(1\) & M0 &                \(4\) &        \(1\) & \(1\) &   \((\alpha_2)\) & \(0\) \\
 \(1\) & \(28\) &  \(84\) & \(3\) &  \(-6\) &       \(1\) &        \(\theta\) &        \(4.431\) & \((96)\)           &    \(1\) & M0 &                \(4\) &        \(4\) & \(3\) & \((\alpha_1^3)\) & \(0\) \\
 \(1\) & \(29\) &  \(87\) & \(9\) &  \(-2\) &       \(0\) &        \(\theta\) &        \(4.466\) & \((30)\)           &    \(1\) & M0 &                \(4\) &        \(1\) & \(1\) &   \((\alpha_2)\) & \(0\) \\
\hline
\end{tabular}
\end{center}
\end{table}


\subsection{Special coefficient \(b=1\) without prime divisors}
\label{ss:One}

\noindent
This is the only situation where a non-trivial absolute principal factor does \textit{not necessarily} occur.

\begin{example}
\label{exm:One}
In Table
\ref{tbl:One},
in fact, only for \(r=1\), an absolute DPF with norm \(3\) exists.
The generator \(\alpha=2+\theta^2\) of the prime ideal
\(\alpha\mathcal{O}_L=\mathfrak{Q}\in\mathbb{P}_L\) with \(3\mathcal{O}_L=\mathfrak{Q}^3\) is sporadic.
The fundamental unit \(\varepsilon_0\) is always the zero \(\theta\) of the polynomial \(P(X)=X^3+3rX-1\)
with norm \(\mathrm{N}_{L/\mathbb{Q}}(\theta)=b=1\).
Always when \(r\equiv -1\,(\mathrm{mod}\,3)\),
the fundamental discriminant is \(d_K\equiv +1\,(\mathrm{mod}\,3)\),
we necessarily have conductor \(f=9\),
the prime \(3\) splits in \(K\),
and \(L\) is of type \(\alpha_2\) with relative DPF.
On the other hand,
if \(r\equiv +1\,(\mathrm{mod}\,3)\),
the fundamental discriminant is \(d_K\equiv +3\,(\mathrm{mod}\,9)\),
we necessarily have conductor \(f=3\),
and \(L\) is of type \(\alpha_1\) with capitulation.
The multiplicity is \(m=2^{\tau-1}=1\) for \(\varrho=0\),
\(m=3^\varrho=3\) for \(\varrho=1\), and
\(m=3^\varrho=9\) for \(\varrho=2\).
\end{example}


\renewcommand{\arraystretch}{1.1}
\begin{table}[ht]
\caption{Conductors \(f\) with a single prime divisor and tiny multiplets}
\label{tbl:OneReal}
\begin{center}
\begin{tabular}{|r|r|r||c|c|c|c|r||r|r|r|c|c|}
\hline
 \(b\) & \(r\)   & \(a\)    & \(f\) & \(d_K\) & \(\varrho\) & \(\varepsilon_1\) & \(\mathrm{Reg}\) & \(m^{\prime\prime}\) & \(m^\prime\) & \(m\) & Multiplet         & \(v\) \\
\hline
 \(1\) &  \(-1\) &   \(-3\) & \(9\) &     --- &         --- &        \(\theta\) &              --- &                  --- &          --- &   --- &       \((\zeta)\) &   --- \\
\hline
 \(1\) &  \(-2\) &   \(-6\) & \(3\) &   \(3\) &       \(0\) &        \(\theta\) &        \(6.801\) &                \(1\) &        \(1\) & \(1\) & \((\varepsilon)\) & \(0\) \\
 \(1\) &  \(-4\) &  \(-12\) & \(9\) &   \(4\) &       \(0\) &   \(\theta^2-12\) &        \(8.951\) &                \(4\) &        \(1\) & \(1\) &    \((\delta_2)\) & \(0\) \\
 \(1\) &  \(-5\) &  \(-15\) & \(3\) &   \(3\) &       \(0\) &   \(\theta^2-15\) &       \(17.936\) &                \(1\) &        \(1\) & \(1\) & \((\varepsilon)\) & \(0\) \\
 \(1\) &  \(-7\) &  \(-21\) & \(9\) &   \(7\) &       \(0\) &        \(\theta\) &       \(12.539\) &                \(4\) &        \(1\) & \(1\) &    \((\delta_2)\) & \(0\) \\
 \(1\) &  \(-8\) &  \(-24\) & \(3\) &   \(3\) &       \(0\) &   \(\theta^2-24\) &       \(64.883\) &                \(1\) &        \(1\) & \(1\) & \((\varepsilon)\) & \(0\) \\
 \(1\) & \(-10\) &  \(-30\) & \(9\) &   \(1\) &       \(0\) &   \(\theta^2-30\) &       \(38.143\) &                \(4\) &        \(1\) & \(1\) &     \((\beta_2)\) & \(0\) \\
 \(1\) & \(-11\) &  \(-33\) & \(3\) &   \(3\) &       \(0\) &        \(\theta\) &       \(50.076\) &                \(1\) &        \(1\) & \(1\) & \((\varepsilon)\) & \(0\) \\
 \(1\) & \(-13\) &  \(-39\) & \(9\) &   \(4\) &       \(0\) &        \(\theta\) &       \(22.915\) &                \(4\) &        \(1\) & \(1\) &    \((\delta_2)\) & \(0\) \\
 \(1\) & \(-16\) &  \(-48\) & \(9\) &   \(7\) &       \(0\) &        \(\theta\) &       \(41.424\) &                \(4\) &        \(1\) & \(1\) &    \((\delta_2)\) & \(0\) \\
 \(1\) & \(-17\) &  \(-51\) & \(3\) &   \(3\) &       \(1\) &   \(\theta^2-51\) &       \(28.523\) &                \(4\) &        \(4\) & \(3\) &  \((\delta_1^3)\) & \(0\) \\
 \(1\) & \(-19\) &  \(-57\) & \(9\) &   \(1\) &       \(0\) &        \(\theta\) &       \(67.400\) &                \(4\) &        \(1\) & \(1\) &     \((\beta_2)\) & \(0\) \\
 \(1\) & \(-20\) &  \(-60\) & \(3\) &   \(3\) &       \(1\) &   \(\theta^2-60\) &      \(121.987\) &                \(4\) &        \(4\) & \(3\) &  \((\delta_1^3)\) & \(0\) \\
 \(1\) & \(-22\) &  \(-66\) & \(9\) &   \(4\) &       \(1\) &   \(\theta^2-66\) &       \(32.042\) &               \(13\) &        \(4\) & \(3\) &  \((\alpha_2^3)\) & \(0\) \\
 \(1\) & \(-23\) &  \(-69\) & \(3\) &   \(3\) &       \(1\) &   \(\theta^2-69\) &       \(38.472\) &                \(4\) &        \(4\) & \(3\) &  \((\delta_1^3)\) & \(0\) \\
 \(1\) & \(-25\) &  \(-75\) & \(9\) &   \(7\) &       \(0\) &        \(\theta\) &       \(24.090\) &                \(4\) &        \(1\) & \(1\) &    \((\delta_2)\) & \(0\) \\
 \(1\) & \(-26\) &  \(-78\) & \(3\) &   \(3\) &       \(0\) &   \(\theta^2-78\) &      \(245.988\) &                \(1\) &        \(1\) & \(1\) & \((\varepsilon)\) & \(0\) \\
 \(1\) & \(-28\) &  \(-84\) & \(9\) &   \(1\) &       \(0\) &        \(\theta\) &      \(182.518\) &                \(4\) &        \(1\) & \(1\) &     \((\beta_2)\) & \(0\) \\
 \(1\) & \(-29\) &  \(-87\) & \(3\) &   \(3\) &       \(1\) &        \(\theta\) &       \(25.415\) &                \(4\) &        \(4\) & \(3\) & \((\beta_1^2,\delta_1)\) & \(0\) \\
\hline
\end{tabular}
\end{center}
\end{table}


\section{Experimental verification of totally real cubic fields}
\label{s:ExperimentsReal}

\noindent
We continue with \textit{totally real cubic} fields, \(\sigma=-1\), in the Sections \S\S\
\ref{ss:OneReal}--\ref{ss:FiveReal}.


\subsection{The special coefficient \(b=1\) without prime divisors}
\label{ss:OneReal}

\noindent
This is the only situation where a non-trivial absolute principal factor does \textit{not necessarily} occur.
\begin{example}
\label{exm:OneReal}
Nevertheless, an absolute DPF with norm \(3\) exists for the types
\(\zeta\), \(\varepsilon\), \(\beta_2\) in Table
\ref{tbl:OneReal}.
The generator \(\alpha\in\mathcal{O}_L\) of the prime ideal
\(\alpha\mathcal{O}_L=\mathfrak{Q}\in\mathbb{P}_L\) with \(3\mathcal{O}_L=\mathfrak{Q}^3\) is sporadic
(not parametrized), for instance
\(\alpha=\theta^2-\theta-3\) for \(r=-1\),
\(\alpha=2+\theta\) for \(r=-2\),
\(\alpha=4-\theta\) for \(r=-5\),
\(\alpha=8+205\theta-42\theta^2\) for \(r=-8\),
\(\alpha=2\theta^2-11\theta\) for \(r=-10\),
\(\alpha=2\theta^2-11\theta-3\) for \(r=-11\),
\(\alpha=\theta^2+9\theta+11\) for \(r=-19\),
\(\alpha=3\,099\,202\theta^2+27\,391\,256\theta+350\,661\) for \(r=-26\),
\(\alpha=6492\theta^2+59\,539\theta+710\) for \(r=-28\).
The columns are arranged in the same way as in Table
\ref{tbl:Five},
except that the unique fundamental unit \(\varepsilon_0\)
is replaced by one of the generators \(\varepsilon_1\)
of a fundamental system \(U_L=\langle -1,\varepsilon_1,\varepsilon_2\rangle\) with two units.
The sign \(\sigma=-1\) is absorbed in \(r<0\).
Exceptionally, the field \(L\) with \(r=-1\) is a unique cyclic cubic field
without associated quadratic field \(K\).
One of the two generating units, \(\varepsilon_1\),
is either the zero \(\theta\) of the polynomial \(P(X)=X^3+3rX-1\)
with norm \(\mathrm{N}_{L/\mathbb{Q}}(\theta)=b=1\)
or the inverse \(\theta^{-1}=\theta^2+a\) (because \(\theta\cdot(\theta^2+a)=\theta^3+a\theta=+1\)).
When \(r\equiv -1\,(\mathrm{mod}\,3)\),
the fundamental discriminant is \(d_K\equiv +1\,(\mathrm{mod}\,3)\),
we necessarily have conductor \(f=9\),
the prime \(3\) splits in \(K\),
and \(L\) is of type \(\beta_2\) or \(\delta_2\) with relative DPF,
in the case of positive \(3\)-class rank \(\varrho=1\) also
\(\alpha_2\) with additional capitulation (for \(r=-22\)).
When \(r\equiv +1\,(\mathrm{mod}\,3)\),
the fundamental discriminant is \(d_K\equiv +3\,(\mathrm{mod}\,9)\),
we necessarily have conductor \(f=3\),
and \(L\) is of type \(\varepsilon\) or \(\delta_1\) with capitulation
when \(\varrho=1\) (but \textit{not} of type \(\beta_1\) for \(r=-29\)).
In both cases, \(\tau=1\), the conductor \(f\) is \textit{regular} with \(\omega=0\),
and \(K\) is real quadratic with unit group \(U_K=\langle -1,\eta\rangle\) of rank one,
fundamental unit \(\eta\), and \(3\)-Selmer rank \(\sigma_3=\varrho_3+1\).
Even for \(\varrho_3=0\), the restrictive Formula
\eqref{eqn:Multiplicity}
must be considered.
However, since \((u,v)=(0,1)\) leads to \(m=0\),
only the \textit{free} situation \((u,v)=(1,0)\) yields a positive multiplicity
\(m=3^{\varrho+\omega}\cdot 2^u\cdot\frac{1}{3}\left(2^{v-1}-(-1)^{v-1}\right)
=3^{0+0}\cdot 2^1\cdot\frac{1}{3}\left(2^{-1}-(-1)^{-1}\right)=1\cdot 2\cdot\frac{1}{2}=1\).
For \(\varrho_3=1\), we directly use
\cite[Thm 1.1, p. 832]{Ma1992},
\(m^\prime=(3^{\varrho+\tau+\omega-\delta}-1)/2=(3^{1+1+0-\delta}-1)/2=m(1)+m(f)\),
which implies \(m^\prime=m(1)=1\), \(m(f)=0\) for \(\delta=1\),
and \(m^\prime=4\), \(m(1)=1\), \(m(f)=3\) for \(\delta=0\) (and thus \(v=0\)).
\end{example}


\subsection{A prime coefficient \(b=q\)}
\label{ss:TwoReal}

\begin{example}
\label{exm:TwoReal}
In Table
\ref{tbl:TwoReal},
we consider \(b=7\) with \textit{regular} conductor \(f=21=3\cdot 7\).
\end{example}


\renewcommand{\arraystretch}{1.1}
\begin{table}[ht]
\caption{Regular conductor \(f\) with \(\tau=2\) and modest multiplets}
\label{tbl:TwoReal}
\begin{center}
\begin{tabular}{|r|r|r||c|c|c|c|r|c||r|r|r|c|c|}
\hline
 \(b\) & \(r\)   & \(a\)    & \(f\)  & \(d_K\) & \(\varrho\) & \(\varepsilon_1\)      & \(\mathrm{Reg}\) & \(\mathrm{Cl}(L)\) & \(m^{\prime\prime}\) & \(m^\prime\) & \(m\) & Multiplet             & \(v\) \\
\hline
 \(7\) &  \(-2\) &  \(-42\) & \(21\) &   \(3\) &       \(0\) &   \((1+6\theta)^{-1}\) &       \(48.890\) & \((3)\)            &                \(4\) &        \(1\) & \(1\) &         \((\beta_2)\) & \(2\) \\
 \(7\) &  \(-3\) &  \(-63\) & \(21\) &   \(6\) &       \(0\) & \(1+8\theta+\theta^2\) &       \(32.497\) & \((6)\)            &               \(13\) &        \(4\) & \(2\) & \((\varepsilon^{2})\) & \(0\) \\
 \(7\) &  \(-5\) & \(-105\) & \(21\) &   \(3\) &       \(0\) &\(44+6\theta-\theta^2\) &       \(67.313\) & \((6)\)            &                \(4\) &        \(1\) & \(1\) &         \((\beta_2)\) & \(2\) \\
\hline
\end{tabular}
\end{center}
\end{table}


\subsection{A composite coefficient \(b=q_1\cdots q_k\)}
\label{ss:FiveReal}

\begin{example}
\label{exm:ThreeReal}
In Table
\ref{tbl:ThreeReal},
we consider \(b=6=3\cdot 2\) with \textit{irregular} conductor \(f=18=3^2\cdot 2\),
\(\omega=1\), \(\tau=2\).
Here, we must pay attention to the degenerate situation \(\delta(3)=1\),
where the multiplicity \(m=3^\varrho\cdot 2^{\tau-1}\)
is independent of \(v\in\lbrace 0,1,2\rbrace\), according to
\cite[Thm. 3.4, p. 2217]{Ma2014}.
Actually, this situation
occurs for \(\varrho=0\), whence \(m=2\).
For \(r\in\lbrace -1,-13,-14\rbrace\), we obviously have \(\delta(3)=0\), and we can use Formula
\eqref{eqn:Multiplicity},
which yields \(m=3^{0+1}\cdot 2^0\cdot 1=3\) with \((u,v)=(0,2)\).
We point out that the \(3\)-ring spaces are given by
\(V(3)=V\), \(V(2)=V(9)=0\), for \(r\in\lbrace -1,-13,-14\rbrace\).
The multiplicity \(m=3\) for \(r\in\lbrace -8,-15\rbrace\) with \(\varrho=1\)
is neither divisible by \(9\) nor by \(2\).
This enforces \(\delta=2\) with \(\delta(3)=1\), \(u=n=0\), \(v=2\) in
\cite[Thm. 4.2, p. 2225]{Ma2014}.
For \(r\in\lbrace -2,-8,-18\rbrace\), the unit \(\varepsilon_1\) is sporadic.
The free situation arises for \(r=-18\) with \(m=3^{0+1}\cdot 2^{2-1}=6\).
\end{example}


\renewcommand{\arraystretch}{1.1}
\begin{table}[ht]
\caption{Irregular conductor \(f\) with \(\tau=2\) and modest multiplets}
\label{tbl:ThreeReal}
\begin{center}
\begin{tabular}{|r|r|r||c|c|c|c|r|c||r|r|r|c|c|}
\hline
 \(b\) & \(r\)   & \(a\)    & \(f\)  & \(d_K\) & \(\varrho\) & \(\varepsilon_1\)      & \(\mathrm{Reg}\) & \(\mathrm{Cl}(L)\) & \(m^{\prime\prime}\) & \(m^\prime\) & \(m\) & Multiplet             & \(v\) \\
\hline
 \(6\) &  \(-1\) &  \(-18\) & \(18\) &   \(6\) &       \(0\) &          \(1+3\theta\) &       \(37.349\) & \(1\)              &               \(13\) &        \(4\) & \(3\) &      \((\gamma^{3})\) & \(2\) \\
 \(6\) &  \(-2\) &  \(-36\) & \(18\) &   \(6\) &       \(0\) & \(1+6\theta+\theta^2\) &       \(17.800\) & \((6)\)            &               \(13\) &        \(4\) & \(2\) & \((\varepsilon^{2})\) & \(0\) \\
 \(6\) &  \(-3\) &  \(-54\) & \(18\) &   \(6\) &       \(0\) &          \(1+9\theta\) &      \(112.244\) & \(1\)              &               \(13\) &        \(4\) & \(2\) &      \((\gamma^{2})\) & \(1\) \\
 \(6\) &  \(-4\) &  \(-72\) & \(18\) &   \(6\) &       \(0\) &  \((1+12\theta)^{-1}\) &      \(439.653\) & \(1\)              &               \(13\) &        \(4\) & \(2\) &      \((\gamma^{2})\) & \(1\) \\
 \(6\) &  \(-5\) &  \(-90\) & \(18\) &   \(6\) &       \(0\) &         \(1+15\theta\) &      \(338.141\) & \(1\)              &               \(13\) &        \(4\) & \(2\) &      \((\gamma^{2})\) & \(1\) \\
 \(6\) &  \(-6\) & \(-108\) & \(18\) &   \(6\) &       \(0\) &         \(1+18\theta\) &      \(120.692\) & \((3)\)            &               \(13\) &        \(4\) & \(2\) & \((\varepsilon^{2})\) & \(0\) \\
 \(6\) &  \(-7\) & \(-126\) & \(18\) &   \(6\) &       \(0\) &         \(1+21\theta\) &     \(1053.567\) & \(1\)              &               \(13\) &        \(4\) & \(2\) &      \((\gamma^{2})\) & \(1\) \\
 \(6\) &  \(-8\) & \(-144\) & \(18\) &   \(6\) &       \(1\) &      \(287-2\theta^2\) &       \(39.919\) & \((24)\)           &               \(13\) &        \(4\) & \(3\) & \((\varepsilon^{3})\) & \(2\) \\
 \(6\) &  \(-9\) & \(-162\) & \(18\) &   \(6\) &       \(0\) &         \(1+27\theta\) &      \(167.251\) & \((6)\)            &               \(13\) &        \(4\) & \(2\) & \((\varepsilon^{2})\) & \(0\) \\
 \(6\) & \(-10\) & \(-180\) & \(18\) &   \(6\) &       \(0\) &         \(1+30\theta\) &      \(403.368\) & \((3)\)            &               \(13\) &        \(4\) & \(2\) &      \((\gamma^{2})\) & \(0\) \\
 \(6\) & \(-13\) & \(-234\) & \(18\) &   \(6\) &       \(0\) &         \(1+39\theta\) &      \(449.354\) & \((2)\)            &               \(13\) &        \(4\) & \(3\) &      \((\gamma^{3})\) & \(2\) \\
 \(6\) & \(-14\) & \(-252\) & \(18\) &   \(6\) &       \(0\) &  \((1+42\theta)^{-1}\) &      \(572.305\) & \((5)\)            &               \(13\) &        \(4\) & \(3\) &      \((\gamma^{3})\) & \(2\) \\
 \(6\) & \(-15\) & \(-270\) & \(18\) &   \(6\) &       \(1\) &  \((1+45\theta)^{-1}\) &      \(679.670\) & \((3)\)            &               \(13\) &        \(4\) & \(3\) &     \((\beta_1^{3})\) & \(2\) \\
 \(6\) & \(-18\) & \(-324\) & \(18\) &   \(6\) &       \(0\) &\(1+54\theta-3\theta^2\)&       \(56.862\) & \((2,24)\)         &               \(40\) &       \(13\) & \(6\) & \((\varepsilon^{6})\) & \(0\) \\
\hline
\end{tabular}
\end{center}
\end{table}


\begin{example}
\label{exm:FourReal}
In Table
\ref{tbl:FourReal},
the field \(L\) is of type \(\beta_2\) with relative DPF for \(r=-3\),
of type \(\beta_1\) with capitulation for \(r=-24\),
but otherwise exclusively of type \(\gamma\) with
\textit{two-dimensional absolute differential principal factorization}
(in particular also for \(r\in\lbrace -11,-17\rbrace\)).
The conductor \(f=210\) is constant,
and the columns are arranged in the same way as in Table
\ref{tbl:Five},
except that \(\varepsilon_0\)
is replaced by \(\varepsilon_1\)
in \(U_L=\langle -1,\varepsilon_1,\varepsilon_2\rangle\).
The norm \(b=70\) of the absolute DPF \(\theta\) of \(L/\mathbb{Q}\) and
one generating unit of \(L\), \(\varepsilon_1=1-3r\theta\) or \((1-3r\theta)^{-1}\),
verify Theorem
\ref{thm:TrinomialPrinciple},
except for \(r=-24\) with
\(\varepsilon_1=1+71\theta+\theta^2\),
\(\varepsilon_2=367919+\theta-73\theta^2\).
The conductor is \textit{regular} with \(\tau=4\) and \(\omega=0\).
In the case \(\varrho=0\), we can use Formula
\eqref{eqn:Multiplicity},
\(m=3^{\varrho+\omega}\cdot 2^u\cdot\frac{1}{3}\left(2^{v-1}-(-1)^{v-1}\right)
=3^{0+0}\cdot 2^u\cdot\frac{1}{3}\left(2^{v-1}-(-1)^{v-1}\right)\),
and we get
\begin{itemize}
\item
\(m=16\cdot\frac{1}{2}=8\) with \((u,v)=(4,0)\) for no value of \(r\) (the conductor is \textit{never free}),
\item
\(m=4\cdot 1=4\) with \((u,v)=(2,2)\) for \(r\in\lbrace -6,-11,-15,-18,-20,-21,-23,-30\rbrace\),
\item
\(m=2\cdot 1=2\) with \((u,v)=(1,3)\) for \(r\in\lbrace -2,-3,-14,-17,-27\rbrace\),
\item
\(m=1\cdot 3=3\) with \((u,v)=(0,4)\) for \(r\in\lbrace -5,-9,-26\rbrace\).
\end{itemize}
In the case  \(\varrho=1\), we first use the accumulative formula
\(m^\prime=(3^{\varrho+\tau+\omega-\delta}-1)/2=(3^{1+4+0-\delta}-1)/2\)
in order to get the defect \(\delta\) from the experiment:
\begin{itemize}
\item
\(m^\prime=(3^{5}-1)/2=121\) with \(\delta=0\) never occurs,
\item
\(m^\prime=(3^{4}-1)/2=40\) with \(\delta=1\) occurs only once for \(r=-24\),
\item
\(m^\prime=(3^{3}-1)/2=13\) with \(\delta=2\) occurs twice for \(r\in\lbrace -12,-29\rbrace\).
\end{itemize}
For \(r=-24\), we can use Formula
\eqref{eqn:Multiplicity},
\(m=3^{\varrho+\omega}\cdot 2^u\cdot\frac{1}{3}\left(2^{v-1}-(-1)^{v-1}\right)
=3^{1+0}\cdot 2^u\cdot\frac{1}{3}\left(2^{v-1}-(-1)^{v-1}\right)=3\cdot 2^0\cdot 3=9\),
which enforces \((u,v)=(0,4)\).
For \(r\in\lbrace -12,-29\rbrace\) we have to employ the formula for defect \(\delta=2\) in
\cite[Thm. 4.1, p. 2224, and Tbl. 7, p. 2225]{Ma2014},
\(m=3^{\varrho+\omega}\cdot 2^u\cdot\frac{1}{9}\left(2^{v-1}+\sum_{i=1}^{4}\,(-1)^{v-n_i}2^{n_i}\right)=6\),
which is not unambiguous, since it is realized for
either \((u,v)=(1,3)\), \((n_1,n_2,n_3,n_4)=(1,1,1,0)\)
or \((u,v)=(0,4)\), \((n_1,n_2,n_3,n_4)=(2,2,0,0)\).
\end{example}


\renewcommand{\arraystretch}{1.0}
\begin{table}[ht]
\caption{Regular conductor \(f\) with \(\tau=4\) and modest multiplets}
\label{tbl:FourReal}
\begin{center}
\begin{tabular}{|r|r|r||c|c|c|c|r||r|r|r|c|c|}
\hline
 \(b\)  & \(r\)   & \(a\)     & \(f\)   & \(d_K\) & \(\varrho\) & \(\varepsilon_1\)     & \(\mathrm{Reg}\) & \(m^{\prime\prime}\) & \(m^\prime\) & \(m\) & Multiplet         & \(v\) \\
\hline
 \(70\) &  \(-2\) &  \(-420\) & \(210\) &   \(3\) &       \(0\) &         \(1+6\theta\) &      \(382.970\) &               \(40\) &       \(13\) & \(2\) &  \((\gamma^{2})\) & \(3\) \\
 \(70\) &  \(-3\) &  \(-630\) & \(210\) &   \(6\) &       \(0\) &  \((1+9\theta)^{-1}\) &      \(239.661\) &               \(40\) &       \(13\) & \(2\) & \((\beta_2^{2})\) & \(3\) \\
 \(70\) &  \(-5\) & \(-1050\) & \(210\) &   \(3\) &       \(0\) & \((1+15\theta)^{-1}\) &     \(2096.691\) &               \(40\) &       \(13\) & \(3\) &  \((\gamma^{3})\) & \(4\) \\
 \(70\) &  \(-6\) & \(-1260\) & \(210\) &   \(6\) &       \(0\) &        \(1+18\theta\) &      \(598.126\) &               \(40\) &       \(13\) & \(4\) &  \((\gamma^{4})\) & \(2\) \\
 \(70\) &  \(-9\) & \(-1890\) & \(210\) &   \(6\) &       \(0\) &        \(1+27\theta\) &     \(2170.666\) &               \(40\) &       \(13\) & \(3\) &  \((\gamma^{3})\) & \(4\) \\
 \(70\) & \(-11\) & \(-2310\) & \(210\) &   \(3\) &       \(0\) &        \(1+33\theta\) &      \(617.937\) &               \(40\) &       \(13\) & \(4\) & \((\beta_2^2,\gamma^2)\) & \(2\) \\
 \(70\) & \(-12\) & \(-2520\) & \(210\) &   \(6\) &       \(1\) & \((1+36\theta)^{-1}\) &     \(6448.722\) &               \(40\) &       \(13\) & \(6\) &  \((\gamma^{6})\) & \(3,4\) \\
 \(70\) & \(-14\) & \(-2940\) & \(210\) &   \(3\) &       \(0\) & \((1+42\theta)^{-1}\) &      \(598.130\) &               \(40\) &       \(13\) & \(2\) &  \((\gamma^{2})\) & \(3\) \\
 \(70\) & \(-15\) & \(-3150\) & \(210\) &   \(6\) &       \(0\) &        \(1+45\theta\) &    \(10084.464\) &               \(40\) &       \(13\) & \(4\) &  \((\gamma^{4})\) & \(2\) \\
 \(70\) & \(-17\) & \(-3570\) & \(210\) &   \(3\) &       \(0\) &        \(1+51\theta\) &     \(3989.244\) &               \(40\) &       \(13\) & \(2\) & \((\beta_2,\gamma)\) & \(3\) \\
 \(70\) & \(-18\) & \(-3780\) & \(210\) &   \(6\) &       \(0\) &        \(1+54\theta\) &     \(2888.134\) &               \(40\) &       \(13\) & \(4\) &  \((\gamma^{4})\) & \(2\) \\
 \(70\) & \(-20\) & \(-4200\) & \(210\) &   \(3\) &       \(0\) & \((1+60\theta)^{-1}\) &     \(7545.737\) &               \(40\) &       \(13\) & \(4\) &  \((\gamma^{4})\) & \(2\) \\
 \(70\) & \(-21\) & \(-4410\) & \(210\) &   \(6\) &       \(0\) & \((1+63\theta)^{-1}\) &    \(18219.191\) &               \(40\) &       \(13\) & \(4\) &  \((\gamma^{4})\) & \(2\) \\
 \(70\) & \(-23\) & \(-4830\) & \(210\) &   \(3\) &       \(0\) &        \(1+69\theta\) &    \(16124.360\) &               \(40\) &       \(13\) & \(4\) &  \((\gamma^{4})\) & \(2\) \\
 \(70\) & \(-24\) & \(-5040\) & \(210\) &   \(6\) &       \(1\) & \(1+71\theta+\theta^2\) &    \(121.037\) &              \(121\) &       \(40\) & \(9\) & \((\beta_1^{9})\) & \(4\) \\
 \(70\) & \(-26\) & \(-5460\) & \(210\) &   \(3\) &       \(0\) &        \(1+78\theta\) &    \(30813.917\) &               \(40\) &       \(13\) & \(3\) &  \((\gamma^{3})\) & \(4\) \\
 \(70\) & \(-27\) & \(-5670\) & \(210\) &   \(6\) &       \(0\) & \((1+81\theta)^{-1}\) &    \(22405.140\) &               \(40\) &       \(13\) & \(2\) &  \((\gamma^{2})\) & \(3\) \\
 \(70\) & \(-29\) & \(-6090\) & \(210\) &   \(3\) &       \(1\) & \((1+87\theta)^{-1}\) &    \(25155.090\) &               \(40\) &       \(13\) & \(6\) &  \((\gamma^{6})\) & \(3,4\) \\
 \(70\) & \(-30\) & \(-6300\) & \(210\) &   \(6\) &       \(0\) & \((1+90\theta)^{-1}\) &    \(13416.106\) &               \(40\) &       \(13\) & \(4\) &  \((\gamma^{4})\) & \(2\) \\
\hline
\end{tabular}
\end{center}
\end{table}


\begin{example}
\label{exm:FiveReal}
In Table
\ref{tbl:FiveReal},
the field \(L\) is exclusively of type \(\gamma\),
in particular also for \(r\in\lbrace -7,-20,-21\rbrace\).
The conductor \(f=6930=3^2\cdot 770\) is \textit{irregular} with \(\omega=1\),
since \(d_K\equiv -3\,(\mathrm{mod}\,9)\),
and \(\tau=5\).
The accumulative formula for the constant value
\(m^\prime=(3^{\varrho+\tau+\omega-\delta}-1)/2=(3^{\varrho+5+1-\delta}-1)/2=121\)
yields the defect \(\delta=\varrho+1\) from the experiment.
Thus, \(\delta=1\) for \(\varrho=0\),
and we can use Formula
\eqref{eqn:Multiplicity}
in the few cases with \(\delta_3(3)=0\),
i.e. when the prime \(3\) is free.
While \((u,v)=(4,1)\) is meanwhile well-known to be impossible, we have the following admissible pairs:
\begin{itemize}
\item
\((u,v)=(3,2)\) with \(m=3^{0+1}\cdot 2^{3}\cdot 1=24\) for \(r=-2\),
\item
\((u,v)=(2,3)\) with \(m=3^{0+1}\cdot 2^{2}\cdot 1=12\) for \(r=-14\),
\item
\((u,v)=(1,4)\) with \(m=3^{0+1}\cdot 2\cdot 3=18\) for \(r=-17\),
\item
\((u,v)=(0,5)\) with \(m=3^{0+1}\cdot 1\cdot 5=15\) for \(r=-20\) (with \textit{five restrictive} primes).
\end{itemize} 
The dominating majority, however, has \(\delta_3(3)=1\),
i.e. the prime \(3\) is restrictive,
and we must use the degenerate formula 
\(m=3^\varrho\cdot 2^{\tau-1}=3^0\cdot 2^4=16\) in
\cite[Thm. 3.4, p. 2217]{Ma2014}.
There remain the four cases with \(\varrho=1\) and \(\delta=2\):
for \(r\in\lbrace -12,-13\rbrace\),
\cite[Thm. 4.1, p. 2224, and Tbl. 7, p. 2225]{Ma2014}
\(m=3^{\varrho+\omega}\cdot 2^u\cdot\frac{1}{9}\left(2^{v-1}+\sum_{i=1}^{4}\,(-1)^{v-n_i}2^{n_i}\right)=9\cdot 2=18\),
with \(\delta_3(3)=0\) can be realized in three ways,
either \((u,v)=(1,4)\), \((n_1,n_2,n_3,n_4)=(2,1,1,0)\)
or \((u,v)=(0,5)\), \((n_1,n_2,n_3,n_4)=(2,1,1,1)\)
or \((u,v)=(0,5)\), \((n_1,n_2,n_3,n_4)=(3,2,0,0)\).
For  \(\delta_3(3)=1\), however, we need the degenerate formula 
\(m=3^{\varrho}\cdot 2^{u_{\mathrm{eff}}}\cdot\frac{1}{3}\left(2^{v_{\mathrm{eff}}-1}-(-1)^{v_{\mathrm{eff}}-1}\right)\)
\cite[Thm. 4.2, p. 2225]{Ma2014}.
This yields \(m=3^1\cdot 2^{2}\cdot 1=12\) with \((u_{\mathrm{eff}},v_{\mathrm{eff}})=(2,3)\) for \(r=-24\), and
\(m=3^1\cdot 2^{3}\cdot 1=24\) with \((u_{\mathrm{eff}},v_{\mathrm{eff}})=(3,2)\) for \(r=-11\),
in the notation of
\cite[Rmk. 4.2, p. 2225]{Ma2014}.
\end{example}


\begin{remark}
\label{rmk:DPFtypes}
The prediction of the possible \textit{types of differential principal factorizations} in our Main Theorem
\ref{thm:Main}
has been fully verified for all cubic number fields \(L=\mathbb{Q}(\theta)\),
generated by a zero \(\theta\) of a monogenic trinomial
\(P(X)=X^3+aX-b\in\mathbb{Z}\lbrack X\rbrack\) with \(b\ge 2\) and \(a=\sigma\cdot 3rb\):
we have exclusively DPF-type \(\beta\) for the simply real fields in the Tables
\ref{tbl:Five}--\ref{tbl:Two},
and we have DPF-types \(\beta_2\) and \(\varepsilon\) in Table
\ref{tbl:TwoReal},
and DPF-types \(\gamma\), \(\varepsilon\) and \(\beta_1\) in Table
\ref{tbl:ThreeReal},
for the totally real fields.
When the number of prime divisors \(\tau\) of the conductor \(f\) increases,
the tendency is clearly towards dominance of DPF-type \(\gamma\).

Several tables have revealed quite frequent occurrence of a precise
\(3\)-class group \(\mathrm{Cl}(L)\simeq\mathbb{Z}/3\mathbb{Z}\) of order \(3\),
which we exploit in applications to homogeneous splitting in the next Section \S\
\ref{s:Application}.
\end{remark}


\renewcommand{\arraystretch}{1.0}
\begin{table}[ht]
\caption{Irregular conductor \(f\) with \(\tau=5\) and most extensive multiplets}
\label{tbl:FiveReal}
\begin{center}
\begin{tabular}{|r|r|r||c|c|c|r|r||r|r|r|c|}
\hline
 \(b\)    & \(r\)   & \(a\)       & \(f\)    & \(d_K\) & \(\varrho\) & \(\varepsilon_1\)     & \(\mathrm{Reg}\) & \(m^{\prime\prime}\) & \(m^\prime\) & \(m\)  & Multiplet         \\
\hline
 \(2310\) &  \(-1\) &   \(-6930\) & \(6930\) &   \(6\) &       \(0\) &         \(1+3\theta\) &     \(3162.092\) &             \(1093\) &      \(121\) & \(16\) & \((\gamma^{16})\) \\
 \(2310\) &  \(-2\) &  \(-13860\) & \(6930\) &   \(6\) &       \(0\) &  \((1+6\theta)^{-1}\) &    \(51307.528\) &             \(1093\) &      \(121\) & \(24\) & \((\gamma^{24})\) \\
 \(2310\) &  \(-3\) &  \(-20790\) & \(6930\) &   \(6\) &       \(0\) &         \(1+9\theta\) &     \(4473.254\) &             \(1093\) &      \(121\) & \(16\) & \((\gamma^{16})\) \\
 \(2310\) &  \(-4\) &  \(-27720\) & \(6930\) &   \(6\) &       \(0\) & \((1+12\theta)^{-1}\) &    \(97292.964\) &             \(1093\) &      \(121\) & \(16\) & \((\gamma^{16})\) \\
 \(2310\) &  \(-5\) &  \(-34650\) & \(6930\) &   \(6\) &       \(0\) &        \(1+15\theta\) &   \(123805.784\) &             \(1093\) &      \(121\) & \(16\) & \((\gamma^{16})\) \\
 \(2310\) &  \(-6\) &  \(-41580\) & \(6930\) &   \(6\) &       \(0\) &        \(1+18\theta\) &   \(205926.744\) &             \(1093\) &      \(121\) & \(16\) & \((\gamma^{16})\) \\
 \(2310\) &  \(-7\) &  \(-48510\) & \(6930\) &   \(6\) &       \(0\) &        \(1+21\theta\) &   \(157999.413\) &             \(1093\) &      \(121\) & \(16\) & \((\beta_2^6,\gamma^{10})\) \\
 \(2310\) &  \(-8\) &  \(-55440\) & \(6930\) &   \(6\) &       \(0\) & \((1+24\theta)^{-1}\) &   \(226072.445\) &             \(1093\) &      \(121\) & \(16\) & \((\gamma^{16})\) \\
 \(2310\) &  \(-9\) &  \(-62370\) & \(6930\) &   \(6\) &       \(0\) &        \(1+27\theta\) &    \(14047.524\) &             \(1093\) &      \(121\) & \(16\) & \((\gamma^{16})\) \\
 \(2310\) & \(-10\) &  \(-69300\) & \(6930\) &   \(6\) &       \(0\) & \((1+30\theta)^{-1}\) &   \(104108.440\) &             \(1093\) &      \(121\) & \(16\) & \((\gamma^{16})\) \\
 \(2310\) & \(-11\) &  \(-76230\) & \(6930\) &   \(6\) &       \(1\) & \((1+33\theta)^{-1}\) &    \(72766.909\) &             \(1093\) &      \(121\) & \(24\) & \((\gamma^{24})\) \\
 \(2310\) & \(-12\) &  \(-83160\) & \(6930\) &   \(6\) &       \(1\) & \((1+36\theta)^{-1}\) &   \(146025.443\) &             \(1093\) &      \(121\) & \(18\) & \((\gamma^{18})\) \\
 \(2310\) & \(-13\) &  \(-90090\) & \(6930\) &   \(6\) &       \(1\) & \((1+39\theta)^{-1}\) &   \(435685.904\) &             \(1093\) &      \(121\) & \(18\) & \((\gamma^{18})\) \\
 \(2310\) & \(-14\) &  \(-97020\) & \(6930\) &   \(6\) &       \(0\) & \((1+42\theta)^{-1}\) &   \(118606.194\) &             \(1093\) &      \(121\) & \(12\) & \((\gamma^{12})\) \\
 \(2310\) & \(-15\) & \(-103950\) & \(6930\) &   \(6\) &       \(0\) & \((1+45\theta)^{-1}\) &   \(788927.927\) &             \(1093\) &      \(121\) & \(16\) & \((\gamma^{16})\) \\
 \(2310\) & \(-16\) & \(-110880\) & \(6930\) &   \(6\) &       \(0\) & \((1+48\theta)^{-1}\) &   \(137853.305\) &             \(1093\) &      \(121\) & \(16\) & \((\gamma^{16})\) \\
 \(2310\) & \(-17\) & \(-117810\) & \(6930\) &   \(6\) &       \(0\) & \((1+51\theta)^{-1}\) &   \(634328.382\) &             \(1093\) &      \(121\) & \(18\) & \((\gamma^{18})\) \\
 \(2310\) & \(-18\) & \(-124740\) & \(6930\) &   \(6\) &       \(0\) &        \(1+54\theta\) &    \(67005.254\) &             \(1093\) &      \(121\) & \(16\) & \((\gamma^{16})\) \\
 \(2310\) & \(-19\) & \(-131670\) & \(6930\) &   \(6\) &       \(0\) & \((1+57\theta)^{-1}\) &   \(281980.092\) &             \(1093\) &      \(121\) & \(16\) & \((\gamma^{16})\) \\
 \(2310\) & \(-20\) & \(-138600\) & \(6930\) &   \(6\) &       \(0\) &        \(1+60\theta\) &  \(1082708.016\) &             \(1093\) &      \(121\) & \(15\) & \((\beta_2^5,\gamma^{10})\) \\
 \(2310\) & \(-21\) & \(-145530\) & \(6930\) &   \(6\) &       \(0\) &        \(1+63\theta\) &   \(319723.332\) &             \(1093\) &      \(121\) & \(16\) & \((\beta_2^{15},\gamma)\) \\
 \(2310\) & \(-22\) & \(-152460\) & \(6930\) &   \(6\) &       \(0\) &        \(1+66\theta\) &   \(446353.410\) &             \(1093\) &      \(121\) & \(16\) & \((\gamma^{16})\) \\
 \(2310\) & \(-23\) & \(-159390\) & \(6930\) &   \(6\) &       \(0\) & \((1+69\theta)^{-1}\) &  \(1060309.933\) &             \(1093\) &      \(121\) & \(16\) & \((\gamma^{16})\) \\
 \(2310\) & \(-24\) & \(-166320\) & \(6930\) &   \(6\) &       \(1\) & \((1+72\theta)^{-1}\) &   \(262045.082\) &             \(1093\) &      \(121\) & \(12\) & \((\gamma^{12})\) \\
\hline 
\end{tabular}
\end{center}
\end{table}

\newpage

\section{Application to homogeneous splitting}
\label{s:Application}

\noindent
In this applied section,
we demonstrate the high practical value of our cubic trinomials
\(P(X)=X^3+aX-b\)
in numerical examples verifying theoretical statements.
Sections \S\S\
\ref{s:Experiments},
\ref{s:ExperimentsReal}
revealed a considerable proportion of cubic fields \(L\)
with non-trivial \(3\)-class group \(\mathrm{Cl}_3(L)\).
Therefore, let \(L\) be a cubic number field
with maximal order \(\mathcal{O}_L\)
and ideal class group \(\mathrm{Cl}(L)\) of exact order \(3\).
Denote the class of an ideal \(\mathfrak{a}\) of \(\mathcal{O}_L\)
by \(\lbrack\mathfrak{a}\rbrack\).
The following concept was coined by S. Chavan
\cite[Dfn. 1.1]{Ch2022}.

\begin{definition}
\label{dfn:HomoSplit}
The field \(L\) is said to possess \textit{homogeneous splitting},
if for every prime number \(p\in\mathbb{P}\) which splits completely in \(L\),
that is \(p\mathcal{O}_L=\mathfrak{p}_1\mathfrak{p}_2\mathfrak{p}_3\),
the three distinct prime ideals \(\mathfrak{p}_1,\mathfrak{p}_2,\mathfrak{p}_3\in\mathbb{P}_L\)
belong to the same class,
\(\lbrack\mathfrak{p}_1\rbrack=\lbrack\mathfrak{p}_2\rbrack=\lbrack\mathfrak{p}_3\rbrack\).
\end{definition}

\noindent
Note that this phenomenon is trivial for a cyclic cubic field \(L\)
with \(\mathrm{Gal}(L/\mathbb{Q})=\langle\sigma\rangle\), since
\(\mathfrak{p}_2=\mathfrak{p}_1^\sigma\),
\(\mathfrak{p}_3=\mathfrak{p}_1^{\sigma^2}\),
and thus
\(\mathfrak{p}_1=\pi\mathcal{O}_L\) implies \(\mathfrak{p}_2=\pi^\sigma\mathcal{O}_L\),
\(\mathfrak{p}_3=\pi^{\sigma^2}\mathcal{O}_L\) and
\(\lbrack\mathfrak{p}_1\rbrack=\lbrack\mathfrak{p}_2\rbrack=\lbrack\mathfrak{p}_3\rbrack=1\).
Similarly for \(\lbrack\mathfrak{p}_1\rbrack\ne 1\),
always using the general relation
\(\lbrack\mathfrak{p}_1\rbrack\cdot\lbrack\mathfrak{p}_2\rbrack\cdot\lbrack\mathfrak{p}_3\rbrack=1\)
in the class group \(\mathrm{Cl}(L)\).

Let \(H\) be the Hilbert class field of \(L\),
that is the maximal abelian unramified extension field of \(L\).
According to the assumption that \(\mathrm{Cl}(L)\simeq\mathbb{Z}/3\mathbb{Z}\),
and by the Artin reciprocity law \(\mathrm{Gal}(H/L)\simeq\mathrm{Cl}(L)\),
the field \(H\) is of absolute degree
\(\lbrack H:\mathbb{Q}\rbrack=\lbrack H:L\rbrack\cdot\lbrack L:\mathbb{Q}\rbrack=3\cdot 3=9\)
over the rational number field \(\mathbb{Q}\).
Denote by \(S\) the normal closure of \(H\).

In
\cite[Thm. 1.2]{Ch2022},
Chavan has proved that
\(S=H\) if and only if \(L\) is a cyclic cubic field, and
\(L\) possesses homogeneous splitting if and only if the relative degree
\(\lbrack S:H\rbrack\) is either \(1\) or \(2\).
With other words:
A non-Galois cubic number field \(L\) has homogeneous splitting if and only if
the absolute degree
\(\lbrack S:\mathbb{Q}\rbrack=\lbrack S:H\rbrack\cdot\lbrack H:\mathbb{Q}\rbrack=2\cdot 9\)
is \(18\).

The sufficiency of this criterion was proved by means of
the decomposition law for prime ideals in the cubic field \(L\)
and the class field prime factorization in \(H\).
The proof of the necessity requires an overview of the possible degrees
\(\lbrack S:H\rbrack\in\lbrace 1,2,6,18\rbrace\),
the classification of transitive permutation groups \(\mathcal{P}\le\mathcal{S}_n\)
of degree \(n=9\) and orders \(\#\mathcal{P}\in\lbrace 9,18,54,162\rbrace\),
and the Chebotarev density theorem
in order to prove the existence of a prime \(p\in\mathbb{P}\)
with non-homogeneous splitting in \(L\) for \(\#\mathcal{P}\ge 54\).

On the first glance,
it is not clear if the rather restrictive constraint of a class group \(\mathrm{Cl}(L)\)
with precise order \(3\) can be realized for infinite families of non-Galois cubic fields \(L\),
and it arouses ones curiosity how the possible orders
\(\#\mathrm{Gal}(S/\mathbb{Q})\in\lbrace 18,54,162\rbrace\)
are actually distributed for non-Galois cubic fields \(L\).

We selected some of our monogenic trinomials of the form
\(P(X)=X^3+aX-b\in\mathbb{Z}\lbrack X\rbrack\) with \(a=3rb\)
such that the maximal order \(\mathcal{O}_L=\mathbb{Z}\lbrack\theta\rbrack\)
of the cubic field \(L=\mathbb{Q}(\theta)\)
has a power basis \((1,\theta,\theta^2)\)
in terms of the traceless real zero \(\theta\) with \(P(\theta)=0\).
Computations were done with Magma
\cite{BCP1997,BCFS2022,MAGMA2022},
based on the class field routines by Fieker
\cite{Fi2001}.

The following tables give the conductor \(f\),
discriminant \(d_L\),
regulator \(\mathrm{Reg}\), and
the absolute degree \(\lbrack S:\mathbb{Q}\rbrack\)
of the normal closure \(S\) of the Hilbert class field \(H\) of \(L\). 

The required condition of a class group \(\mathrm{Cl}(L)\simeq\mathbb{Z}/3\mathbb{Z}\)
is obviously limited to finitely many cases of simply-real cubic number fields \(L\)
with negative discriminant \(d_L<0\), e.g. those in Table
\ref{tbl:ExamplesSimplyReal}.

\renewcommand{\arraystretch}{1.1}
\begin{table}[ht]
\caption{Complete splitting in simply-real cubic number fields \(L\)}
\label{tbl:ExamplesSimplyReal}
\begin{center}
\begin{tabular}{|r||r|r|r||r|r|r|c|c|}
\hline
 No.   & \(b\) & \(r\)  & \(a\)  & \(f\)  & \(d_L\)    & \(\mathrm{Reg}\) & \(\lbrack S:\mathbb{Q}\rbrack\) & Splitting \\
\hline
 \(1\) & \(1\) &  \(2\) &  \(6\) &  \(9\) &   \(-891\) &        \(1.796\) &                         \(18\)  & homogeneous \\
 \(2\) & \(1\) &  \(5\) & \(15\) &  \(9\) & \(-13527\) &        \(2.708\) &                         \(18\)  & homogeneous \\
 \(3\) & \(2\) &  \(2\) & \(12\) &  \(6\) &  \(-7020\) &        \(6.075\) &                         \(54\)  & non-homogeneous \\
 \(4\) & \(2\) &  \(3\) & \(18\) &  \(6\) & \(-23446\) &        \(7.287\) &                         \(54\)  & non-homogeneous \\
 \(5\) & \(5\) &  \(1\) & \(15\) & \(45\) & \(-14175\) &        \(4.927\) &                         \(54\)  & non-homogeneous \\
 \(6\) & \(7\) &  \(1\) & \(21\) & \(21\) & \(-38367\) &        \(5.257\) &                         \(18\)  & homogeneous \\
\hline
\end{tabular}
\end{center}
\end{table}

\noindent
For No. \(1\), the primes \(p\in\lbrace 23,89\rbrace\) split completely in the form
\(p\mathcal{O}_L=\mathfrak{p}_1\mathfrak{p}_2\mathfrak{p}_3\) in \(L\),
and the splitting is homogeneous,
since for \(p=23\), all prime ideals \(\mathfrak{p}_i\in\mathbb{P}_L\)
lie in the same non-principal class \(\lbrack\mathfrak{p}_i\rbrack\ne 1\),
whereas for \(p=89\), all prime ideals belong to the principal class \(\lbrack\mathfrak{p}_i\rbrack=1\).

For No. \(3\), the primes \(p\in\lbrace 11,17,61,97\rbrace\) split completely in the form
\(p\mathcal{O}_L=\mathfrak{p}_1\mathfrak{p}_2\mathfrak{p}_3\) in \(L\),
but the splitting is non-homogeneous, since for each \(p\),
one of the prime ideals \(\mathfrak{p}_i\in\mathbb{P}_L\) belongs to the principal class
and the other two prime ideals lie in distinct non-principal classes.

The required condition of a class group \(\mathrm{Cl}(L)\simeq\mathbb{Z}/3\mathbb{Z}\)
seems to occur for infinitely many cases of totally-real cubic number fields \(L\)
with positive discriminant \(d_L>0\), e.g. those in Table
\ref{tbl:ExamplesTotallyReal}.

\renewcommand{\arraystretch}{1.1}
\begin{table}[ht]
\caption{Complete splitting in totally-real cubic number fields \(L\)}
\label{tbl:ExamplesTotallyReal}
\begin{center}
\begin{tabular}{|r||r|r|r||r|r|r|c|c|}
\hline
 No.   & \(b\) & \(r\)   & \(a\)   & \(f\)  & \(d_L\)    & \(\mathrm{Reg}\) & \(\lbrack S:\mathbb{Q}\rbrack\) & Splitting \\
\hline
 \(1\) & \(1\) &  \(-4\) & \(-12\) &  \(9\) &   \(6885\) &        \(8.951\) &                         \(18\)  & homogeneous \\
 \(2\) & \(1\) &  \(-7\) & \(-21\) &  \(9\) &  \(37017\) &       \(12.539\) &                         \(18\)  & homogeneous \\
 \(3\) & \(1\) & \(-10\) & \(-30\) &  \(9\) & \(107973\) &       \(38.143\) &                         \(18\)  & homogeneous \\
 \(4\) & \(1\) & \(-13\) & \(-39\) &  \(9\) & \(237249\) &       \(22.915\) &                         \(18\)  & homogeneous \\
 \(5\) & \(1\) & \(-16\) & \(-48\) &  \(9\) & \(442341\) &       \(41.424\) &                         \(18\)  & homogeneous \\
 \(6\) & \(1\) & \(-17\) & \(-51\) &  \(3\) & \(530577\) &       \(28.523\) &                         \(54\)  & non-homogeneous \\
 \(7\) & \(1\) & \(-19\) & \(-57\) &  \(9\) & \(740745\) &       \(67.400\) &                         \(18\)  & homogeneous \\
 \(8\) & \(1\) & \(-20\) & \(-60\) &  \(3\) & \(863973\) &      \(121.987\) &                         \(54\)  & non-homogeneous \\
\hline
\end{tabular}
\end{center}
\end{table}

\noindent
For No. \(7\), the primes \(p\in\lbrace 11,19,47,53,97\rbrace\) split completely in the form
\(p\mathcal{O}_L=\mathfrak{p}_1\mathfrak{p}_2\mathfrak{p}_3\) in \(L\),
and the splitting is homogeneous,
since for \(p\in\lbrace 11,47,97\rbrace\), all prime ideals \(\mathfrak{p}_i\in\mathbb{P}_L\)
lie in the same non-principal class \(\lbrack\mathfrak{p}_i\rbrack\ne 1\),
whereas for \(p\in\lbrace 19,53\rbrace\), they belong to the principal class \(\lbrack\mathfrak{p}_i\rbrack=1\).

For No. \(8\), the primes \(p\in\lbrace 19,29,97\rbrace\) split completely in the form
\(p\mathcal{O}_L=\mathfrak{p}_1\mathfrak{p}_2\mathfrak{p}_3\) in \(L\),
but the splitting is non-homogeneous, since for each \(p\),
one of the prime ideals \(\mathfrak{p}_i\in\mathbb{P}_L\) belongs to the principal class
and the other two prime ideals lie in distinct non-principal classes.

We could not find examples with \(\lbrack S:\mathbb{Q}\rbrack=162\).
Either their occurrence is extremely sparse,
and an explicit example by the author would be illuminating,
or they cannot occur for theoretical reasons.
This seems to be an open problem for future clarification.


\section{Conclusion}
\label{s:Conclusion}

\noindent
In this article,
we have shown that a cubic number field \(L=\mathbb{Q}(\theta)\)
generated by a zero \(\theta\) of a \textit{monogenic} trinomial
\(X^3+aX-b\in\mathbb{Z}\lbrack X\rbrack\)
with squarefree absolute coefficient \(b\in\mathbb{N}\) and
\(a=\sigma\cdot 3rb\), \(\sigma\in\lbrace -1,+1\rbrace\), \(r\in\mathbb{N}\),
is \textit{distinguished} by special arithmetical properties
within its multiplet \((L_1,\ldots,L_m)\) which consists of all
non-isomorphic cubic fields sharing a common conductor \(f=3^e\cdot b\)
and a common associated quadratic fundamental discriminant \(d_K\)
such that \(d_L=d_{L_i}=f^2\cdot d_K\) for all \(1\le i\le m\).
The exceptional number theoretic invariants of \(L\) are
a \textit{parametrized absolute differential principal factorization} (DPF)
\(\mathfrak{Q}=\theta\mathcal{O}_L\) such that \(b\mathcal{O}_L=\mathfrak{Q}^3\), and,
in the case of a \textit{simply real} field with \(\sigma=+1\),
a \textit{parametrized fundamental unit} \(\varepsilon_0=1-3r\theta\)
such that \(U_L=\langle -1,\varepsilon_0\rangle\).
It is particularly remarkable,
that in the case of positive \(3\)-class rank \(\varrho_K\ge 1\),
which enables \textit{capitulation} of \(K=\mathbb{Q}(\sqrt{d_L})\) in \(N=L\cdot K\),
and in the case of a prime divisor \(q\equiv +1\,(\mathrm{mod}\,3)\) of \(b\),
splitting in \(K\),
which admits a \textit{relative} DPF of \(N/K\),
the other two options (capitulation and relative DPF)
are rigorously supressed by the a priori absolute DPF.

The \textit{monogeneity} of \(L\),
providing a \textit{power integral basis} of \(\mathcal{O}_L\),
was a decisive tool in many proofs.
The crucial break through, however, is due to our revival of
the Ukrainian mathematician Georgi F. Voronoi,
born \(1886\) in Zhuravka near Kiev,
who coined the concept of \textit{lattice minima}
and designed his marvellous \textit{Voronoi algorithm} for the
construction of \textit{chains} of lattice minima.
Without his techniques,
we would neither have been able to prove the fundamentality of
the parametrized unit \(\varepsilon_0=1-3r\theta\)
nor the striking explicit chain with \textit{very small period length}
\(\ell\in\lbrace 1,2,3\rbrace\)
and the sharp lower bound \(r\ge L(b)=\lceil b/3\rceil\)
for \(\ell=3\) with an \(\mathrm{M}2\)-chain.

All experimental results in \S\S\
\ref{s:Experiments},
\ref{s:ExperimentsReal}
are in perfect accordance with our theoretical predictions.
We conclude with a conjecture concerning the
splitting field \(N\) of the trinomial \(P(X)=X^p+aX-b\),
which generalizes Corollary
\ref{cor:Symmetric},
when \(p\ge 5\).

\begin{conjecture}
\label{cnj:Symmetric}
For any odd prime number \(p\ge 5\),
the splitting field \(N\) of the trinomial \(P(X)=X^p+aX-b\)
with \(a=\sigma\cdot prb\), \(\sigma\in\lbrace -1,+1\rbrace\),
\(r\in\mathbb{N}\), and squarefree \(b\in\mathbb{N}\)
has the full symmetric automorphism group
\(\mathrm{Gal}(N/\mathbb{Q})\simeq S_p\).
\end{conjecture}


\section{Acknowledgements}
\label{s:Gratifications}

\noindent
The first author acknowledges that his research was supported by
the Austrian Science Fund (FWF): projects P26008-N25 and J0497-PHY,
and by the Research Executive Agency of the European Union (EUREA).



\end{document}